\documentclass[twocolumn]{autart}    

\usepackage{graphicx}          

\usepackage{tabularx}
\usepackage{natbib}
\usepackage{colortbl}
\usepackage{booktabs}
\usepackage{algorithm}
\usepackage{algpseudocode}

\usepackage{amsmath,amssymb,enumerate,subfigure,multirow,epsfig,psfrag,hhline,color,srcltx,subeqnarray}
\usepackage{multirow}
\usepackage{slashbox}
\usepackage{graphicx}          

\usepackage{color}
\usepackage{xcolor}
\definecolor{dblue}{rgb}{0,0,7}
\definecolor{dred}{rgb}{7,0,0}
\definecolor{dgreen}{rgb}{0,2,0}
\definecolor{mycolor}{rgb}{0.75,0,0.6}
\usepackage{hyperref}
\hypersetup{breaklinks=true,colorlinks=true,citecolor=blue}

\def\cal{\mathcal}

\newcolumntype{Z}{>{\small\centering\arraybackslash}X}

\newtheorem{theorem}{Theorem}
\newtheorem{corollary}{Corollary}

\newtheorem{lemma}{Lemma}
\newtheorem{example}{Example}

\newtheorem{remark}{Remark}

\allowdisplaybreaks

\begin{document}

\begin{frontmatter}

\title{Periodically time-varying memory static output feedback control design for discrete-time LTI systems} 

\author[Yonsei]{Dong Hwan Lee\thanksref{footnoteinfo2}}\ead{hope2010@yonsei.ac.kr},
\author[Joo]{Young Hoon Joo}\ead{yhjoo@kunsan.ac.kr},
\author[Joo]{Myung Hwan Tak}\ead{takgom@kunsan.ac.kr}
\thanks[footnoteinfo2]{Corresponding author. Tel.: +82 2 2123
2773; fax: +82 2 362 4539.}
\address[Yonsei]{Department of Electrical and Computer Engineering, Purdue University, West Lafayette, IN 47906, USA}
\address[Joo]{Department of Control and Robotics Engineering, Kunsan National University, Kunsan, Chonbuk, 573-701,
Korea}


\begin{abstract}                          
This paper addresses the problem of static output feedback (SOF)
stabilization for discrete-time LTI systems. We approach this
problem using the recently developed periodically time-varying
memory state-feedback controller (PTVMSFC) design scheme. A
bilinear matrix inequality (BMI) condition which uses a
pre-designed PTVMSFC is developed to design the periodically
time-varying memory SOF controller (PTVMSOFC). The BMI condition
can be solved by using BMI solvers. Alternatively, we can apply
two-steps and iterative linear matrix inequality algorithms that
alternate between the PTVMSFC and PTVMSOFC designs. Finally, an
example is given to illustrate the proposed methods.
\end{abstract}

\end{frontmatter}

\section{Introduction}

The design of static output feedback (SOF) controllers has
received a significant amount of attention to date since it is
common experience in practical control applications that having
full access to the state is not always possible. While a wide
variety of problems related to controller analysis and design can
be recast as convex linear matrix inequality (LMI) problems
\citep*{Boyd1994} which are easily tractable by standard convex
optimization techniques \citep*{Gahinet1995, Strum1999,
Lofberg2004}, this is not the case for the SOF problem
\citep*{Fu1997} since the most general characterization of the SOF
design is bilinear matrix inequalities (BMIs) for which complete
and efficient methods to find their global solutions are not
available yet. For this reason, the SOF design is one of the most
challenging open problems in the control literature. Nowadays,
there is immense literature addressing the SOF problem through
various approaches, just to name a few:
\begin{itemize}
\item A simple method using structural properties of the open-loop
system \citep*{Garcia2001};

\item Iterative schemes based on the linear quadratic regulator
(LQR) theory \citep*{Kucera1995,Rosinova2003};

\item Sufficient LMI conditions using similarly transformations
\citep*{Prempain2001,Lee2006, Dong2007} and using the elimination
lemma \citep*{Dong2013};

\item Sufficient LMI conditions with linear matrix equality
constraints \citep*{Crusius1999};

\item Two-steps LMI approaches \citep*{Bara2005} using a
congruence transformation and fixing the Lyapunov matrix
structure;

\item Iterative LMI (ILMI) methods based on the LQR theory
\citep*{Cao1998}, cone complementarity linearization
\citep*{Ghaoui1997,He2006}, quadratic separation concept
\citep*{Peaucelle2005}, descriptor system augmentation
\citep*{Shu2010}, and substitutive ILMI algorithm
\citep*{Fujimori2004};

\item ILMI schemes \citep*{Peaucelle2001} and two-steps LMI
approaches \citep*{Mehdi2004,Agulhari2010, Agulhari2012}
alternating between state-feedback (SF) and SOF designs;

\item Mixed LMI/randomized methods \citep*{Arzelier2010};

\item The rank constrained LMI strategy \citep*{Orsi2006};

\item Nonlinear optimization approaches
\citep*{Goh1994,Kanev2004,Henrion2005,Burke2006}.
\end{itemize}

In this paper, we consider the problem of designing a SOF
controller for discrete-time LTI systems. Among the important
results mentioned earlier, the main idea of this paper is
motivated by \cite{Peaucelle2001}, where efficient ILMI procedures
that alternate between the SF and SOF designs are developed based
on the elimination lemma \citep{Boyd1994}. The idea was further
developed in \cite{Mehdi2004} for discrete-time LTI systems by
introducing new decision variables, in \cite{Arzelier2010} in
combination with hit-and-run strategies, and recently in
\cite{Agulhari2010, Agulhari2012} for reduced-order robust ${\cal
H}_\infty$ control of continuous-time uncertain LTI systems.

We revisit this idea in a somewhat different direction for
discrete-time LTI systems. More specifically, our method is an
extension of the work presented in \cite{Peaucelle2001,
Mehdi2004,Agulhari2010, Agulhari2012} to the so-called
periodically time-varying memory controller technique, which was
developed recently by \cite*{Ebihara2009,Ebihara2011,Tregouet2011,
Tregouet2012,Tregouet2013} for robust control purposes.

In the field of robust control of LTI systems, the development of
less conservative robust SF control design has been a fundamental
and challenging problem. In the late 1990s, the so-called extended
Schur complement and slack variable approaches were developed by
the pioneering work in \cite*{de_Oliveira1999, Peaucelle2000,
de_Oliveira2002}, which paved the way for the subsequent
development of the LMI-based robust analysis and control design
approaches (see, e.g, \cite*{Oliveira2007, Oliveira2008} and
references therein). Recently, a new paradigm emerged through a
sequence of interesing researches in \cite{Ebihara2009,
Ebihara2011,Tregouet2011,Tregouet2012,Tregouet2013}, where the
so-called periodically time-varying memory SF controller (PTVMSFC)
which makes use of the state information in a periodic manner was
proposed and turned out to be effective in reducing the
conservatism in the traditional robust SF approaches for
discrete-time systems subject to parameter uncertainties. Despite
those recent progresses, up to the authors' knowledge, an
extension of the PTVMSFC approach to the SOF problem still remains
unresolved.

This paper suggests strategies to design a periodically
time-varying memory SOF controller (PTVMSOFC) that stabilizes
discrete-time LTI systems. To this end, first, we pay attention
for introducing some definitions and notation, which reduce the
difficulty of the matrix calculations and their formal
expressions. Next, by means of the Finsler's lemma
\citep*{Skelton1998}, a necessary and sufficient condition for
designing the PTVMSOFC is derived in terms of BMI problems. Then,
following the lines in \cite{Peaucelle2001}, we use the
elimination lemma \citep*{Skelton1998} to reduce the structure of
the multiplier matrix introduced by the Finsler's lemma to a
special form based on a chosen SF controller, and the BMI problem
comes down to solving another BMI. These BMI problems can be
treated with PENBMI \citep*{Kocvara2005}, a solver for BMIs.
Alternatively, at the price of some conservatism, the BMI problem
can reduce to an LMI problem, based on which the PTVMSOFC design
problem can be solved by applying two-steps LMI and iterative LMI
(ILMI) algorithms \citep{Mehdi2004,Agulhari2010,
Agulhari2012,Peaucelle2001}. Finally, an comparison analysis is
given to evaluate the effectiveness of the proposed approaches.

\section{Preliminaries}

\subsection{Notation}

The adopted notation is as follows: ${\mathbb N}$ and ${\mathbb
N}_+$: sets of nonnegative and positive integers, respectively;
${\mathbb Z}_{[k_1 ,\,k_2 ]}$: set of integers $\{ k_1 ,\,k_1  +
1, \ldots ,\,k_2 \}  \subseteq {\mathbb N}$; ${\mathbb R}^n $:
$n$-dimensional Euclidean space; ${\mathbb R}^{n \times m}$: set
of all $n \times m$ real matrices; ${\mathbb S}_+^n $: set of all
$n \times n$ real symmetric positive definite matrices; $A^T$:
transpose of matrix $A$; ${\rm{He}}\{ A\} : = A^T  + A$; $\rho
(A)$: spectral radius of matrix $A$; $A_ \bot$: any matrices whose
columns form bases of the right null-space of matrix $A$; $A
\otimes B$: Kronecker's product of matrices $A$ and $B$; $A \succ
0$ ($A \prec 0$, $A \succeq 0$, and $A \preceq 0$, respectively):
symmetric positive definite (negative definite, positive
semi-definite, and negative semi-definite, respectively) matrix
$A$; ${\bf{0}}$: zero matrix of appropriate dimensions;
${\bf{0}}_{n \times m} $ and ${\bf{0}}_n $: zero matrix and zero
vector of dimensions $n \times m$ and $n$, respectively;
${\bf{I}}_n $: $n \times n$ identity matrix; ${\mathcal L}_N : =
[\begin{array}{*{20}c}
   {{\bf{I}}_N } & {{\bf{0}}_N }  \\
\end{array}] \in {\mathbb R}^{N \times (N + 1)} $; ${\mathcal R}_N : = [\begin{array}{*{20}c}
   {{\bf{0}}_N} & {{\bf{I}}_N }  \\
\end{array}] \in {\mathbb R}^{N \times (N + 1)} $; ${\bf{e}}_{(N,\,i)} $: unit vector of
dimension $N$ with a $1$ in the $i$-th component and $0$'s
elsewhere; for given two integers $k$ and $N$, $\left\lceil k
\right\rceil_N $: remainder of $k$ divided by $N$;
\begin{align*}
{\bf{T}}_N : = \left[ {\begin{array}{*{20}c}
   0 &  \cdots  & 0 & 1  \\
    \vdots  &  {\mathinner{\mkern2mu\raise1pt\hbox{.}\mkern2mu
 \raise4pt\hbox{.}\mkern2mu\raise7pt\hbox{.}\mkern1mu}}  &  {\mathinner{\mkern2mu\raise1pt\hbox{.}\mkern2mu
 \raise4pt\hbox{.}\mkern2mu\raise7pt\hbox{.}\mkern1mu}}  & 0  \\
   0 & 1 &  {\mathinner{\mkern2mu\raise1pt\hbox{.}\mkern2mu
 \raise4pt\hbox{.}\mkern2mu\raise7pt\hbox{.}\mkern1mu}}  &  \vdots   \\
   1 & 0 &  \cdots  & 0  \\
\end{array}} \right] \in {\mathbb R}^{N \times N}.
\end{align*}

\subsection{Problem formulation}
Consider the discrete-time LTI system described by
\begin{align}
&\left\{ \begin{array}{l}
 x(k + 1) = Ax(k) + Bu(k) \\
 y(k) = Cx(k) \\
 \end{array} \right.\label{system}
\end{align}
where $k \in {\mathbb N}$; $x(k) \in {\mathbb R}^n $ is the state;
$u(k) \in {\mathbb R}^m $ is the control input; $y(k) \in {\mathbb
R}^p $ is the measured output; $\Sigma := (A,\,B,\,C) \in {\mathbb
R}^{n \times n} \times {\mathbb R}^{n \times m}  \times {\mathbb
R}^{p \times n} $ is a tuple of constant matrices. Inspired by the
recently developed PTVMSFC \citep{Ebihara2009,Ebihara2011}, we
suggest the PTVMSOFC (or $N$-PTVMSOFC) of the following form:
\begin{align}
&u(k) = \sum\limits_{i = 0}^{\left\lceil k \right\rceil _N }
{F_{{\rm{SOF}}}^{(\left\lceil k \right\rceil _N ,\,i)} y(k -
i)},\label{PTVMSOFC}
\end{align}
where $N \in {\mathbb N}_+$ is the period of the controller and
$F_{{\rm{SOF}}}^{(\left\lceil k \right\rceil _N ,\,i)}  \in
{\mathbb R}^{m \times p} ,\,(\left\lceil k \right\rceil _N ,\,i)
\in {\mathbb Z}_{[0,\,N - 1]} \times {\mathbb Z}_{[0,\,N - 1]}$
are the SOF gains to be designed. In the case $N=1$, this is the
classical SOF controller. Substituting \eqref{PTVMSOFC} into
\eqref{system}, the $N$-periodic control system (closed-loop
system) can be written as
\begin{align}
&x(k + 1) = Ax(k) + B\sum\limits_{i = 0}^{\left\lceil k
\right\rceil _N } {F_{{\rm{SOF}}}^{(\left\lceil k \right\rceil _N
,\,i)} Cx(k - i)}. \label{SOF-closed-loop-system}
\end{align}
The problem addressed in this paper is to seek the $N$-PTVMSOFC
\eqref{PTVMSOFC} such that the $N$-periodic control system
\eqref{SOF-closed-loop-system} is asymptotically stable.

\section{Main result}
To streamline notation, for two integers $k_1 ,\,k_2  \in {\mathbb
N},\,k_1 \le k_2 $, $x(k_1 :k_2 )$ and $x(k_2 :k_1 )$,
respectively, denote the vectors $x(k_2 :k_1 )^T : =
[\begin{array}{*{20}c}
   {x(k_2 )^T } & {x(k_2  - 1)^T } &  \cdots  & {x(k_1 )^T }  \\
\end{array}]$ and $x(k_1 :k_2 )^T : = [\begin{array}{*{20}c}
   {x(k_1 )^T } & {x(k_1  + 1)^T } &  \cdots  & {x(k_2 )^T }  \\
\end{array}]$.

\subsection{Augmented system representation}

As stated in \cite{Ebihara2011}, for any $k \in \{ k \in {\mathbb
N}:\,\left\lceil k \right\rceil _N = 0\}$, the input of the
PTVMSOFC \eqref{PTVMSOFC} can be expressed in the augmented form:
\begin{align*}
u(k + N - 1:k) =& {\mathcal F}_{{\rm{SOF}}}^{(N,\, \uparrow )} y(k + N - 1:k)\\
=& {\mathcal F}_{{\rm{SOF}}}^{(N,\, \uparrow )} ({\bf I}_N \otimes
C)x(k + N - 1:k),\\
&\forall k \in \{ k \in {\mathbb N}:\,\left\lceil k \right\rceil
_N = 0\},
\end{align*}
where
\begin{align*}
&{\mathcal F}_{{\rm{SOF}}}^{(N,\, \uparrow )} : = \left[
{\begin{array}{*{20}c}
   {F_{{\rm{SOF}}}^{(N - 1,\,0)} } & {F_{{\rm{SOF}}}^{(N - 1,\,1)} } &  \cdots  & {F_{{\rm{SOF}}}^{(N - 1,\,N - 1)} }  \\
   {\bf{0}} &  \ddots  &  \ddots  &  \vdots   \\
    \vdots  &  \ddots  & {F_{{\rm{SOF}}}^{(1,\,0)} } & {F_{{\rm{SOF}}}^{(1,\,1)} }  \\
   {\bf{0}} &  \cdots  & {\bf{0}} & {F_{{\rm{SOF}}}^{(0,\,0)} }  \\
\end{array}} \right].
\end{align*}
In light of this, $N$-periodic control system
\eqref{SOF-closed-loop-system} can be formulated as
\begin{align}
&x(k + N:k + 1) = {\mathcal A}_{{\rm{AUG}}}^{(N,\,\uparrow)} x(k +
N - 1:k),\nonumber\\
&\forall k \in \{ k \in {\mathbb N}:\,\left\lceil k \right\rceil
_N = 0\} , \label{aug-system1}
\end{align}
where ${\mathcal A}_{{\rm{AUG}}}^{(N,\, \uparrow )} : =
({\bf{I}}_N \otimes A) + ({\bf{I}}_N  \otimes B){\mathcal
F}_{{\rm{SOF}}}^{(N,\, \uparrow )} ({\bf{I}}_N \otimes C)$.
Alternatively, based on the transformation $x(k:k + N - 1) =
({\bf{T}}_N  \otimes {\bf{I}}_n )x(k + N - 1:k)$,
\eqref{aug-system1} can be expressed as
\begin{align}
&x(k + 1:k + N) = {\mathcal A}_{{\rm{AUG}}}^{(N,\, \downarrow )}
x(k:k + N - 1),\nonumber\\
&\forall k \in \{ k \in {\mathbb N}:\,\left\lceil k \right\rceil
_N = 0\} ,\label{aug-system2}
\end{align}
where ${\mathcal A}_{{\rm{AUG}}}^{(N,\, \downarrow )} : =
({\bf{I}}_N \otimes A) + ({\bf{I}}_N  \otimes B){\mathcal
F}_{{\rm{SOF}}}^{(N,\, \downarrow )} ({\bf{I}}_N \otimes C)$ and
gain matrix ${\mathcal F}_{{\rm{SOF}}}^{(N,\, \downarrow )}$ takes
the form
\begin{align}
{\mathcal F}_{{\rm{SOF}}}^{(N,\, \downarrow )} :=& ({\bf{T}}_N
\otimes {\bf{I}}_m ){\mathcal F}_{{\rm{SOF}}}^{(N,\, \uparrow )}
({\bf{T}}_N
\otimes {\bf{I}}_p )\nonumber\\
 =& \left[ {\begin{array}{*{20}c}
   {F_{{\rm{SOF}}}^{(0,\,0)} } & {\bf{0}} &  \cdots  & {\bf{0}}  \\
   {F_{{\rm{SOF}}}^{(1,\,1)} } & {F_{{\rm{SOF}}}^{(1,\,0)} } &  \ddots  &  \vdots   \\
    \vdots  &  \ddots  &  \ddots  & {\bf{0}}  \\
   {F_{{\rm{SOF}}}^{(N - 1,\,N - 1)} } &  \cdots  & {F_{{\rm{SOF}}}^{(N - 1,\,1)} } & {F_{{\rm{SOF}}}^{(N - 1,\,0)} }  \\
\end{array}} \right].\label{Fsof}
\end{align}

\subsection{LTI system representation}

Although the underlying system \eqref{system} is an LTI system,
the closed-loop system \eqref{SOF-closed-loop-system} can be
viewed as a periodically time-varying system owing to the
$N$-periodic time-varying controller \eqref{PTVMSOFC}. However, in
order to apply some standard results of LTI systems, it is worth
considering an equivalent LTI representation of
\eqref{SOF-closed-loop-system}. According to \cite{Ebihara2011},
it is always possible that the $N$-periodic control system
\eqref{SOF-closed-loop-system} corresponding to feedback gain
${\cal F}_{{\rm{SOF}}}^{(N,\, \downarrow )}$ can be reformulated
as the equivalent LTI system:
\begin{align}
&\phi (t + 1) = {\cal A}_{{\rm{LTI}}} ({\cal
F}_{{\rm{SOF}}}^{(N,\, \downarrow )} ,\,\Sigma )\phi (t),\quad t
\in {\mathbb N},\label{SOF-LTI-system}
\end{align}
with the state variables $\phi (t) = x(Nt)$, where $\Sigma : =
(A,\,B,\,C)$. Based on the augmented system representation
\eqref{aug-system2}, let
\begin{align*}
&{\cal F}_{{\rm{SOF}}}^{(\delta ,\, \downarrow )} : =
[\begin{array}{*{20}c}
   {{\bf{I}}_{\delta m} } & {{\bf{0}}_{\delta m \times (N - \delta )m} }  \\
\end{array}]{\cal F}_{{\rm{SOF}}}^{(N,\, \downarrow )} [\begin{array}{*{20}c}
   {{\bf{I}}_{\delta p} } & {{\bf{0}}_{\delta p \times (N - \delta )p} }  \\
\end{array}]^T
\end{align*}
and
\begin{align*}
&{\cal A}_{{\rm{AUG}}}^{(\delta ,\, \downarrow )} : =
{\bf{I}}_\delta   \otimes A + ({\bf{I}}_\delta   \otimes B){\cal
F}_{{\rm{SOF}}}^{(\delta ,\, \downarrow )} ({\bf{I}}_\delta
\otimes C),\quad \forall \delta  \in {\mathbb Z}_{[1,\,N]}.
\end{align*}
Furthermore, define ${\cal A}_{{\rm{LTI}}} ({\cal
F}_{{\rm{SOF}}}^{(i,\, \downarrow )} ,\,\Sigma ),\,i \in {\mathbb
Z}_{[1,\,N - 1]}$ as matrices satisfying
\begin{align*}
&x(k + i) = {\cal A}_{{\rm{LTI}}} ({\cal F}_{{\rm{SOF}}}^{(i,\,
\downarrow )} ,\,\Sigma )x(k),\\
&\forall (i,\,k) \in {\mathbb Z}_{[1,\,N - 1]}  \times \{ k \in
{\mathbb N}:\,\left\lceil k \right\rceil _N  = 0\}.
\end{align*}
Then, taking into account \eqref{aug-system2}, it is
straightforward to see that
\begin{align}
x(k + 1:k + N) =&{\cal A}_{{\rm{AUG}}}^{(N,\, \downarrow )}
x(k:k + N - 1)\nonumber\\
=&{\cal A}_{{\rm{AUG}}}^{(N,\, \downarrow )} \left[
{\begin{array}{*{20}c}
   {{\bf{I}}_n }  \\
   {{\cal A}_{{\rm{LTI}}} ({\cal F}_{{\rm{SOF}}}^{(1,\, \downarrow )} ,\,\Sigma )}  \\
    \vdots   \\
   {{\cal A}_{{\rm{LTI}}} ({\cal F}_{{\rm{SOF}}}^{(N - 1,\, \downarrow )} ,\,\Sigma )}  \\
\end{array}} \right]x(k)\nonumber\\
=& \left[ {\begin{array}{*{20}c}
   {{\cal A}_{{\rm{LTI}}} ({\cal F}_{{\rm{SOF}}}^{(1,\, \downarrow )} ,\,\Sigma )}  \\
   {{\cal A}_{{\rm{LTI}}} ({\cal F}_{{\rm{SOF}}}^{(2,\, \downarrow )} ,\,\Sigma )}  \\
    \vdots   \\
   {{\cal A}_{{\rm{LTI}}} ({\cal F}_{{\rm{SOF}}}^{(N,\, \downarrow )} ,\,\Sigma )}  \\
\end{array}} \right]x(k)\nonumber\\
&\forall k \in \{ k \in {\mathbb N}:\,\left\lceil k \right\rceil
_N = 0\},\label{eq1}
\end{align}
and hence, we can obtain the following expression of $ {\cal
A}_{{\rm{LTI}}} ({\cal F}_{{\rm{SOF}}}^{(N,\, \downarrow )}
,\,\Sigma )$:
\begin{align*}
&{\cal A}_{{\rm{LTI}}} ({\cal F}_{{\rm{SOF}}}^{(N,\, \downarrow )}
,\,\Sigma )\\
&\quad = [\begin{array}{*{20}c}
   {{\bf{0}}_{n \times (N - 1)n} } & {{\bf{I}}_n }  \\
\end{array}]{\cal A}_{{\rm{AUG}}}^{(N,\, \downarrow )} \left[ {\begin{array}{*{20}c}
   {{\bf{I}}_n }  \\
   {{\cal A}_{{\rm{LTI}}} ({\cal F}_{{\rm{SOF}}}^{(1,\, \downarrow )} ,\,\Sigma )}  \\
    \vdots   \\
   {{\cal A}_{{\rm{LTI}}} ({\cal F}_{{\rm{SOF}}}^{(N - 1,\, \downarrow )} ,\,\Sigma )}  \\
\end{array}} \right],
\end{align*}
or equivalently, from \eqref{aug-system1}, we have
\begin{align*}
&{\cal A}_{{\rm{LTI}}} ({\cal F}_{{\rm{SOF}}}^{(N,\, \uparrow )}
,\,\Sigma ) = {\cal A}_{{\rm{LTI}}} ({\cal F}_{{\rm{SOF}}}^{(N,\,
\downarrow )} ,\,\Sigma )\\
&\quad = [\begin{array}{*{20}c}
   {{\bf{I}}_n } & {{\bf{0}}_{n \times (N - 1)n} }  \\
\end{array}]{\cal A}_{{\rm{AUG}}}^{(N,\, \uparrow )} \left[ {\begin{array}{*{20}c}
   {{\cal A}_{{\rm{LTI}}} ({\cal F}_{{\rm{SOF}}}^{(N - 1,\, \uparrow )} ,\,\Sigma )}  \\
    \vdots   \\
   {{\cal A}_{{\rm{LTI}}} ({\cal F}_{{\rm{SOF}}}^{(1,\, \uparrow )} ,\,\Sigma )}  \\
   {{\bf{I}}_n }  \\
\end{array}} \right].
\end{align*}
Based on the observation, ${\cal A}_{{\rm{LTI}}} ({\cal
F}_{{\rm{SOF}}}^{(N,\, \downarrow )} ,\,\Sigma )$ can be
constructed using the recursion in Algorithm 1.
\begin{algorithm}[h]
\caption{Construct ${\cal A}_{{\rm{LTI}}} ({\cal
F}_{{\rm{SOF}}}^{(N,\, \downarrow )} ,\,\Sigma )$}
\begin{algorithmic}[1]
\State $\Phi  \leftarrow {\bf{I}}_n$

\For{$\delta  \leftarrow \{ 1,\,2, \ldots ,\,N\}$}

\State ${\cal A}_{{\rm{LTI}}} ({\cal F}_{{\rm{SOF}}}^{(\delta ,\,
\downarrow )} ,\,\Sigma ) \leftarrow [\begin{array}{*{20}c}
   {{\bf{0}}_{n \times (\delta  - 1)n} } & {{\bf{I}}_n }  \\
\end{array}]{\cal A}_{{\rm{AUG}}}^{(\delta ,\, \downarrow )} \Phi$

\State $\Phi  \leftarrow \left[ {\begin{array}{*{20}c}
   \Phi   \\
   {{\cal A}_{{\rm{LTI}}} ({\cal F}_{{\rm{SOF}}}^{(\delta ,\, \downarrow )} ,\,\Sigma )}  \\
\end{array}} \right]$

\EndFor

\State \textbf{return} ${\cal A}_{{\rm{LTI}}} ({\cal
F}_{{\rm{SOF}}}^{(N,\, \downarrow )} ,\,\Sigma )$

\end{algorithmic}
\end{algorithm}

\subsection{$N$-PTVMSOFC synthesis}

We start with the following necessary and sufficient BMI condition
so that \eqref{SOF-LTI-system} is asymptotically stable:
\begin{theorem}\label{theorem1}
There exists ${\cal F}_{{\rm{SOF}}}^{(N,\, \downarrow )} \in
{\mathbb R}^{Nm \times Np}$ defined in \eqref{Fsof} such that
$N$-periodic control system \eqref{SOF-closed-loop-system} or
equivalent LTI system \eqref{SOF-LTI-system} is asymptotically
stable if and only if there exists $P = P^T \in {\mathbb R}^n $
and $M \in {\mathbb R}^{N(n + m) \times ((N - 1)n + Nm)}$ such
that the following problem is satisfied with ${\cal
F}_{{\rm{SOF}}}^{(N,\, \downarrow )} \in {\mathbb R}^{Nm \times
Np}$:
\begin{align}
& P \succ 0,\label{theorem1-eq1}\\
&\Pi _N^T {\cal X}_{N} (P,\,1)\Pi _N  + {\rm{He}}\{ M{\cal
C}({\cal F}_{{\rm{SOF}}}^{(N,\, \downarrow )} )\}  \prec
0,\label{theorem1-eq2}
\end{align}
where
\begin{align*}
&\left\{ \begin{array}{l}
 {\cal X}_\delta  (P,\,\gamma ): = \left( \begin{array}{l}
  - \gamma {\bf{e}}_{(N + 1,\,1)} {\bf{e}}_{(N + 1,\,1)}^T  \\
  + {\bf{e}}_{(N + 1,\,\delta  + 1)} {\bf{e}}_{(N + 1,\,\delta  + 1)}^T  \\
 \end{array} \right) \otimes P; \\
 \Pi _N : = \left[ {\begin{array}{*{20}c}
   {{\bf{e}}_{(N,\,1)}^T  \otimes {\bf{I}}_n } & {{\bf{0}}_{n \times Nm} }  \\
   {{\bf{I}}_N  \otimes A} & {{\bf{I}}_N  \otimes B}  \\
\end{array}} \right]; \\
 {\cal C}({\cal F}_{{\rm{SOF}}}^{(N,\, \downarrow )} ) \\
 : = \left[ {\begin{array}{*{20}c}
   {{\cal F}_{{\rm{SOF}}}^{(N,\, \downarrow )} ({\bf{I}}_N  \otimes C)} & { - {\bf{I}}_N  \otimes {\bf{I}}_m }  \\
   {{\cal L}_{N - 1}  \otimes A - {\cal R}_{N - 1}  \otimes {\bf{I}}_n } & {{\cal L}_{N - 1}  \otimes B}  \\
\end{array}} \right]. \\
 \end{array} \right.
\end{align*}
\end{theorem}
{\bf Proof.} By the Lyapunov argument, $N$-periodic SOF control
system \eqref{SOF-closed-loop-system} is asymptotically stable if
and only if there exists $P \in {\mathbb R}^{n\times n}$ such that
\eqref{theorem1-eq1} and ${\mathcal A}_{{\rm{LTI}}} ({\mathcal
F}_{{\rm{SOF}}}^{(N,\, \downarrow )} )^T P{\mathcal
A}_{{\rm{LTI}}} ({\mathcal F}_{{\rm{SOF}}}^{(N,\, \downarrow )} )
- P \prec 0$ hold. After some algebraic manipulations and using
the relation \eqref{eq1}, one can prove that
\begin{align}
&{\cal A}_{{\rm{LTI}}} ({\cal F}_{{\rm{SOF}}}^{(N,\, \downarrow )}
,\,\Sigma )^T P{\cal A}_{{\rm{LTI}}} ({\cal F}_{{\rm{SOF}}}^{(N,\,
\downarrow )} ,\,\Sigma ) - P\nonumber\\
&= \left[ {\begin{array}{*{20}c}
   {{\bf{I}}_n }  \\
   {{\cal A}_{{\rm{LTI}}} ({\cal F}_{{\rm{SOF}}}^{(1,\, \downarrow )} ,\,\Sigma )}  \\
    \vdots   \\
   {{\cal A}_{{\rm{LTI}}} ({\cal F}_{{\rm{SOF}}}^{(N,\, \downarrow )} ,\,\Sigma )}  \\
\end{array}} \right]^T {\cal X}_N (P,\,1)\left[ {\begin{array}{*{20}c}
   {{\bf{I}}_n }  \\
   {{\cal A}_{{\rm{LTI}}} ({\cal F}_{{\rm{SOF}}}^{(1,\, \downarrow )} ,\,\Sigma )}  \\
    \vdots   \\
   {{\cal A}_{{\rm{LTI}}} ({\cal F}_{{\rm{SOF}}}^{(N,\, \downarrow )} ,\,\Sigma )}  \\
\end{array}} \right]\nonumber\\
&= {\mathcal Q}_N^T \Pi _N^T {\mathcal X}_N (P,\,1)\Pi _N
{\mathcal Q}_N\label{Theorem1-eq3}
\end{align}
and ${\mathcal C}({\mathcal F}_{{\rm{SOF}}}^{(N,\, \downarrow )}
){\mathcal Q}_N  = {\bf{0}}_{((N - 1)n + Nm) \times n}$, where
\begin{align*}
&{\cal Q}_N : = \left[ {\begin{array}{*{20}c}
   {({\bf{I}}_N  \otimes {\bf{I}}_n )}  \\
   {{\cal F}_{{\rm{SOF}}}^{(N,\, \downarrow )} ({\bf{I}}_N  \otimes C)}  \\
\end{array}} \right]\left[ {\begin{array}{*{20}c}
   {{\bf{I}}_n }  \\
   {{\cal A}_{{\rm{LTI}}} ({\cal F}_{{\rm{SOF}}}^{(1,\, \downarrow )} ,\,\Sigma )}  \\
    \vdots   \\
   {{\cal A}_{{\rm{LTI}}} ({\cal F}_{{\rm{SOF}}}^{(N - 1,\, \downarrow )} ,\,\Sigma )}  \\
\end{array}} \right].
\end{align*}
Now, note that ${\mathcal Q}_N$ has full column rank, and
${\rm{rank}}({\mathcal Q}_N ) = n$. Moreover, to use the Finsler's
lemma \citep{Skelton1998}, we need to show that ${\mathcal
C}({\mathcal F}_{{\rm{SOF}}}^{(N,\, \downarrow )} )$ is of full
row rank. To see this, multiplying ${\mathcal C}({\mathcal
F}_{{\rm{SOF}}}^{(N,\, \downarrow )} )$ by the nonsingular matrix
\begin{align*}
&\left[ {\begin{array}{*{20}c}
   {{\cal L}_{N - 1}  \otimes B} & {{\bf{I}}_{(N - 1)n} }  \\
   {{\bf{I}}_{Nm} } & {{\bf{0}}_{Nm \times (N - 1)n} }  \\
\end{array}} \right]
\end{align*}
on the left yields
\begin{align*}
&\left[ {\begin{array}{*{20}c}
   {\left( \begin{array}{l}
 ({\cal L}_{N - 1}  \otimes B){\cal F}_{{\rm{SOF}}}^{(N,\, \downarrow )} ({\bf{I}}_N  \otimes C) \\
  + {\cal L}_{N - 1}  \otimes A - {\cal R}_{N - 1}  \otimes {\bf{I}}_n  \\
 \end{array} \right)} & {\bf{0}}  \\
   {{\cal F}_{{\rm{SOF}}}^{(N,\, \downarrow )} ({\bf{I}}_N  \otimes C)} & { - {\bf{I}}_N  \otimes {\bf{I}}_m }  \\
\end{array}} \right],
\end{align*}
which clearly has full row rank. Therefore,
${\rm{rank}}\,{\mathcal C}({\mathcal F}_{{\rm{SOF}}}^{(N,\,
\downarrow )} ) = (N - 1)n + Nm$ and ${\mathcal C}({\mathcal
F}_{{\rm{SOF}}}^{(N,\, \downarrow )} )$ has a right null-space of
dimension $n$. This implies that ${\mathcal C}({\mathcal
F}_{{\rm{SOF}}}^{(N,\, \downarrow )} )_ \bot = {\mathcal Q}_N$,
and it follows from \eqref{Theorem1-eq3} and
\begin{align*}
&{\cal A}_{{\rm{LTI}}} ({\cal F}_{{\rm{SOF}}}^{(N,\, \downarrow )}
,\,\Sigma )^T P{\cal A}_{{\rm{LTI}}} ({\cal F}_{{\rm{SOF}}}^{(N,\,
\downarrow )} ,\,\Sigma ) - P \prec 0
\end{align*}
that
\begin{align}
&{\cal A}_{{\rm{LTI}}} ({\cal F}_{{\rm{SOF}}}^{(N,\, \downarrow )}
,\,\Sigma )^T P{\cal A}_{{\rm{LTI}}} ({\cal F}_{{\rm{SOF}}}^{(N,\,
\downarrow )} ,\,\Sigma ) - P\nonumber\\
&\quad = {\cal C}({\cal F}_{{\rm{SOF}}}^{(N,\, \downarrow )} )_
\bot ^T \Pi _N^T {\cal X}_N (P,\,1)\Pi _N {\cal C}({\cal
F}_{{\rm{SOF}}}^{(N,\, \downarrow )} )_ \bot\nonumber\\
&\quad \prec 0.\label{Theorem1-eq4}
\end{align}
Applying the Finsler's lemma to \eqref{Theorem1-eq4}, we have that
\eqref{Theorem1-eq4} holds if and only if there exists $M$ such
that \eqref{theorem1-eq2} is satisfied. This completes the proof.
$\blacksquare$

If ${\cal F}_{{\rm{SOF}}}^{(N,\, \downarrow )}$ should be
determined by Theorem \ref{theorem1}, due to the product of
multiplier $M$ introduced by the Finsler's lemma  and controller
parameter ${\mathcal F}_{{\rm{SOF}}}^{(N,\, \downarrow )}$,
\eqref{theorem1-eq2} is a BMI problem. There are several iterative
algorithms to obtain a local solution to BMI problems; for
instance, the alternating minimization algorithm \citep{Goh1994}
is one of the simplest methods. The BMI problem can be also solved
locally by using the BMI solver, PENBMI \citep{Kocvara2005}. It is
important to note that the quality of their solutions depends on
initial parameters of non-convex variables. Therefore, a
reasonable initial guess of the solution can improve the results.
In this context, a very promising result was presented in
\cite{Peaucelle2001}, where based on the {\it a priori} selection
of a suitable SF controller and using elimination lemma
\citep{Boyd1994}, a necessary condition for $M$ to satisfy
\eqref{theorem1-eq2} was derived, and based on this, ILMI
algorithms that alternate between the SF and SOF designs were
proposed. Inspired by the idea in \cite{Peaucelle2001}, we will
suggest an alternative BMI problem which can be viewed as an
extension of those in
\cite{Peaucelle2001,Mehdi2004,Arzelier2010,Agulhari2010,Agulhari2012}.
To this end, we need some preliminary results on the following
PTVMSFC (or $N$-PTVMSFC) proposed in
\cite{Ebihara2009,Ebihara2011}:
\begin{align}
&u(k) = \sum\limits_{i = 0}^{\left\lceil k \right\rceil _N }
{F_{{\rm{SF}}}^{(\left\lceil k \right\rceil _N ,\,i)} x(k -
i)},\label{PTVMSFC}
\end{align}
where $N \in {\mathbb N}_+$ is the period of the controller and
$F_{{\rm{SF}}}^{(\left\lceil k \right\rceil _N ,\,i)} \in {\mathbb
R}^{m \times n} ,\,(\left\lceil k \right\rceil _N ,\,i) \in
{\mathbb Z}_{[0,\,N - 1]} \times {\mathbb Z}_{[0,\,N - 1]}$ are
the SF gains to be designed. Similarly to
\eqref{SOF-closed-loop-system}, substituting \eqref{PTVMSFC} into
\eqref{system} leads to the $N$-periodic SF control system:
\begin{align}
&x(k + 1) = Ax(k) + B\sum\limits_{i = 0}^{\left\lceil k
\right\rceil _N } {F_{{\rm{SF}}}^{(\left\lceil k \right\rceil _N
,\,i)} x(k - i)}. \label{SF-closed-loop-system}
\end{align}
Following the same line as in the PTVMSOFC case, let
\begin{align}
&\xi (t + 1) = {\cal A}_{{\rm{LTI}}} ({\cal F}_{{\rm{SF}}}^{(N,\,
\downarrow )} ,\,\Sigma )\xi (t),\quad t \in {\mathbb
N},\label{SF-LTI-system}
\end{align}
be the equivalent LTI representation of the $N$-periodic SF
control system \eqref{SF-closed-loop-system} corresponding to
PTVMSFC gain matrix ${\mathcal F}_{{\rm{SF}}}^{(N,\, \downarrow
)}$. Then, by using a descriptor-like form of \eqref{system}, the
system-theoretic concept of duality \citep{Ebihara2011}, and the
Finsler's lemma \citep{Skelton1998}, a necessary and sufficient
LMI condition to design \eqref{PTVMSFC} was established in
\cite{Ebihara2011}. For the sake of completeness, it is presented
below.
\begin{lemma}[\cite{Ebihara2011}]\label{Ebihara2011-lemma}
$N$-periodic control system \eqref{SF-closed-loop-system} or
equivalent LTI system \eqref{SF-LTI-system} is asymptotically
stable if and only if there exists matrices $P = P^T  \in {\mathbb
R}^{n\times n} $, $G^{(i,\,j)}  \in {\mathbb R}^{n \times n}$, and
$J^{(i,\,j)}  \in {\mathbb R}^{m \times n}$ such that the
following LMI problem is satisfied:
\begin{align}
&{\cal X}_N (P,\,1) + {\rm{He}}\left\{ \begin{array}{l}
 ({\cal L}_N^T  \otimes A){\cal G} ({\cal R}_N  \otimes {\bf{I}}_n ) \\
  + ({\cal L}_N^T  \otimes B){\cal J} ({\cal R}_N  \otimes {\bf{I}}_n ) \\
  - ({\cal R}_N^T  \otimes {\bf{I}}_n ){\cal G} ({\cal R}_N  \otimes {\bf{I}}_n ) \\
 \end{array} \right\} \prec 0,\label{Ebihara2011-lemma-eq1}
\end{align}
where
\begin{align*}
&{\cal G}: = \left[ {\begin{array}{*{20}c}
   {G^{(1,\,1)} } &  \cdots  & {G^{(1,\,N)} }  \\
   {\bf{0}} &  \ddots  &  \vdots   \\
   {\bf{0}} &  \ddots  & {G^{(N,\,N)} }  \\
\end{array}} \right] \in {\mathbb R}^{Nn \times Nn},\\
&{\cal J}: = \left[ {\begin{array}{*{20}c}
   {J^{(1,\,1)} } &  \cdots  & {J^{(1,\,N)} }  \\
   {\bf{0}} &  \ddots  &  \vdots   \\
   {\bf{0}} &  \ddots  & {J^{(N,\,N)} }  \\
\end{array}} \right] \in {\mathbb R}^{Nm \times Nn}.
\end{align*}
Moreover, an admissible $N$-PTVMSFC gain matrix is given by
\begin{align*}
{\mathcal F}_{{\rm{SF}}}^{(N,\, \uparrow )} : =& \left[
{\begin{array}{*{20}c}
   {F_{{\rm{SF}}}^{(N - 1,\,0)} } & {F_{{\rm{SF}}}^{(N - 1,\,1)} } &  \cdots  & {F_{{\rm{SF}}}^{(N - 1,\,N - 1)} }  \\
   {\bf{0}} &  \ddots  &  \ddots  &  \vdots   \\
    \vdots  &  \ddots  & {F_{{\rm{SF}}}^{(1,\,0)} } & {F_{{\rm{SF}}}^{(1,\,1)} }  \\
   {\bf{0}} &  \cdots  & {\bf{0}} & {F_{{\rm{SF}}}^{(0,\,0)} }  \\
\end{array}} \right]\\
=& {\mathcal J} {\mathcal G}^{-1}.
\end{align*}
\end{lemma}
\begin{remark}\label{theorem1-remark}
Let $(\hat P,\,\hat {\cal J} ,\,\hat {\cal G})$ be a feasible
solution to \eqref{Ebihara2011-lemma-eq1}. Then, $\hat P \succ 0$
holds since pre- and post-multiplying the left-hand side of
\eqref{Ebihara2011-lemma-eq1} by  ${\bf{e}}_{(N + 1,\,1)}^T
\otimes {\bf{I}}_n $ and ${\bf{e}}_{(N + 1,\,1)} \otimes
{\bf{I}}_n $, respectively, results in $P \succ 0$. In addition,
since Lemma \ref{Ebihara2011-lemma} was derived based on a dual
system representation of \eqref{SF-closed-loop-system}, it
guarantees $\hat P \succ 0,\,{\cal A}_{{\rm{LTI}}} (\hat {\cal
F}_{{\rm{SF}}}^{(N,\, \uparrow )} ,\,\Sigma ) \hat P{\cal
A}_{{\rm{LTI}}} (\hat {\cal F}_{{\rm{SF}}}^{(N,\, \uparrow )}
,\,\Sigma )^T - \hat P \prec 0$, and equivalently, $\hat P^{ - 1}
\succ 0$ and
\begin{align*}
&{\cal A}_{{\rm{LTI}}} (\hat {\cal F}_{{\rm{SF}}}^{(N,\,
\downarrow )} ,\,\Sigma )^T \hat P^{ - 1} {\cal A}_{{\rm{LTI}}}
(\hat {\cal F}_{{\rm{SF}}}^{(N,\, \downarrow )} ,\,\Sigma ) - \hat
P^{ - 1}  \prec 0,
\end{align*}
where $\hat {\cal F}_{{\rm{SF}}}^{(N,\, \downarrow )}  =
({\bf{T}}_N  \otimes {\bf{I}}_m )\hat {\cal F}_{{\rm{SF}}}^{(N,\,
\uparrow )} ({\bf{T}}_N  \otimes {\bf{I}}_n )$.
\end{remark}

Before proceeding further, we need to list some definitions. For
any asymptotically stable system matrix $A \in {\mathbb R}^{n
\times n} $, define ${\cal P}(A): = \{ P \in {\mathbb S}_+^n
:\,A^T PA - P \prec 0\} $ as the corresponding set of all
admissible Lyapunov matrices. In addition, let us define
${\mathcal L}_{{\rm{SF}}}^{(N)} (\Sigma ): = \{ {\mathcal
F}_{{\rm{SF}}}^{(N,\, \downarrow )}  \in {\mathbb R}^{Nm \times
Nn} :\,\rho ({\mathcal A}_{{\rm{LTI}}} ({\mathcal
F}_{{\rm{SF}}}^{(N,\, \downarrow )} ,\,\Sigma )) < 1\} $ as the
set of all admissible stabilizing $N$-PTVMSFC gain matrices
${\mathcal F}_{{\rm{SF}}}^{(N,\, \downarrow )}$ corresponding to
$\Sigma $ and ${\mathcal L}_{{\rm{SOF}}}^{(N)} (\Sigma ): = \{
{\mathcal F}_{{\rm{SOF}}}^{(N,\, \downarrow )} \in {\mathbb R}^{Nm
\times Np} :\,\rho ({\mathcal A}_{{\rm{LTI}}} ({\mathcal
F}_{{\rm{SOF}}}^{(N,\, \downarrow )} ,\,\Sigma )\,){\rm{ < 1}}\} $
as the set of all admissible stabilizing $N$-PTVMSOFC gain
matrices corresponding to $\Sigma $. Obviously, ${\mathcal
P}({\mathcal A}_{{\rm{LTI}}} ({\mathcal F}_{{\rm{SF}}}^{(N,\,
\downarrow )} ,\,\Sigma ))$ and ${\mathcal P}({\mathcal
A}_{{\rm{LTI}}} ({\mathcal F}_{{\rm{SOF}}}^{(N,\, \downarrow )}
,\,\Sigma ))$ are the sets of all admissible Lyapunov matrices
corresponding to ${\mathcal F}_{{\rm{SF}}}^{(N,\, \downarrow )}
\in {\mathcal L}_{{\rm{SF}}}^{(N)}$ and ${\mathcal
F}_{{\rm{SOF}}}^{(N,\, \downarrow )}  \in {\mathcal
L}_{{\rm{SOF}}}^{(N)} $, respectively. Now, given a triple $\Sigma
$, let us define
\begin{align*}
&{\cal S}_{{\rm{SF}}}^{(N)} (\Sigma ): = \bigcup\limits_{{\cal
F}_{{\rm{SF}}}^{(N,\, \downarrow )}  \in {\cal
L}_{{\rm{SF}}}^{(N)} (\Sigma )} {{\mathcal P}({\cal
A}_{{\rm{LTI}}}
({\cal F}_{{\rm{SF}}}^{(N,\, \downarrow )} ,\,\Sigma ))},\\
&{\mathcal S}_{{\rm{SOF}}}^{(N)} (\Sigma ): =
\bigcup\limits_{{\mathcal F}_{{\rm{SOF}}}^{(N,\, \downarrow )} \in
{\mathcal L}_{{\rm{SOF}}}^{(N)} (\Sigma )} {{\mathcal P}({\mathcal
A}_{{\rm{LTI}}} ({\mathcal F}_{{\rm{SOF}}}^{(N,\, \downarrow )}
,\,\Sigma ))},
\end{align*}
as the sets of all Lyapunov matrices that corresponds to all
stabilizing PTVMSFC gains ${\mathcal F}_{{\rm{SF}}}^{(N,\,
\downarrow )}  \in {\mathcal L}_{{\rm{SF}}}^{(N)} $ and PTVMSOFC
gains ${\mathcal F}_{{\rm{SOF}}}^{(N,\, \downarrow )}  \in
{\mathcal L}_{{\rm{SOF}}}^{(N)} $, respectively.

Based on the definitions and using the idea that stems from
\cite{Peaucelle2001}, we establish the following theorem:
\begin{theorem}\label{theorem2}
Suppose that $(\hat P,\,\hat {\cal J},\,\hat {\cal G})$ is a
solution to \eqref{Ebihara2011-lemma-eq1}, and let $\hat {\cal
F}_{{\rm{SF}}}^{(N,\, \downarrow )}  = ({\bf{T}}_N \otimes
{\bf{I}}_m )\hat {\cal J}\hat {\cal G}^{- 1} ({\bf{T}}_N \otimes
{\bf{I}}_n )$. Then, there exists matrices ${\cal V} \in {\mathbb
R}^{((N - 1)n + Nm) \times ((N - 1)n + Nm)}$ and ${\cal
F}_{{\rm{SOF}}}^{(N,\, \downarrow )}  \in {\mathbb R}^{Nm \times
Np}$ defined in \eqref{Fsof} such that BMIs \eqref{theorem1-eq1}
and \eqref{theorem1-eq2} in Theorem \ref{theorem1} with $M = {\cal
H}(\hat {\cal F}_{{\rm{SF}}}^{(N,\, \downarrow )} )^T {\cal V}$
have a solution if and only if
\begin{align}
&{\cal P}({\cal A}_{{\rm{LTI}}} (\hat {\cal F}_{{\rm{SF}}}^{(N,\,
\downarrow )} ,\,\Sigma )) \cap {\cal S}_{{\rm{SOF}}}^{(N)}
(\Sigma ) \ne \emptyset\label{nonempty-condition}
\end{align}
holds, where
\begin{align*}
&{\cal H}(\hat {\cal F}_{{\rm{SF}}}^{(N,\, \downarrow )} ): =
\left[ {\begin{array}{*{20}c}
   {\hat {\cal F}_{{\rm{SF}}}^{(N,\, \downarrow )} } & { - {\bf{I}}_N  \otimes {\bf{I}}_m }  \\
   {\left( \begin{array}{l}
 {\cal L}_{N - 1}  \otimes A \\
  - {\cal R}_{N - 1}  \otimes {\bf{I}}_n  \\
 \end{array} \right)} & {{\cal L}_{N - 1}  \otimes B}  \\
\end{array}} \right].
\end{align*}
\end{theorem}
{\bf Proof.} (Sufficiency) If \eqref{nonempty-condition} holds,
then there exists a pair $(P,\,{\cal F}_{{\rm{SOF}}}^{(N,\,
\downarrow )} )$ such that $P \in {\mathbb S}_+^n$,
\begin{align*}
&{\cal A}_{{\rm{LTI}}} ({\cal F}_{{\rm{SOF}}}^{(N,\, \downarrow )}
,\,\Sigma )^T P{\cal A}_{{\rm{LTI}}} ({\cal F}_{{\rm{SOF}}}^{(N,\,
\downarrow )} ,\,\Sigma ) - P \prec 0
\end{align*}
and ${\cal A}_{{\rm{LTI}}} (\hat {\cal F}_{{\rm{SF}}}^{(N,\,
\downarrow )} ,\,\Sigma )^T P{\cal A}_{{\rm{LTI}}} (\hat {\cal
F}_{{\rm{SF}}}^{(N,\, \downarrow )} ,\,\Sigma ) - P \prec 0$ hold.
Following similar lines to the proof of Theorem \ref{theorem1}, we
have that
\begin{align}
&{\cal A}_{{\rm{LTI}}} ({\cal F}_{{\rm{SOF}}}^{(N,\, \downarrow )}
,\,\Sigma )^T P{\cal A}_{{\rm{LTI}}} ({\cal F}_{{\rm{SOF}}}^{(N,\,
\downarrow )} ,\,\Sigma ) - P\nonumber\\
&\quad = {\cal C}({\cal F}_{{\rm{SOF}}}^{(N,\, \downarrow )} )_
\bot ^T \Pi _N^T {\cal X}_N (P,\,1)\Pi _N {\cal C}({\cal
F}_{{\rm{SOF}}}^{(N,\, \downarrow )} )_\bot\nonumber\\
&\quad \prec 0,\label{theorem2-eq1}\\
&{\cal A}_{{\rm{LTI}}} (\hat {\cal F}_{{\rm{SF}}}^{(N,\,
\downarrow )} ,\,\Sigma )^T P{\cal A}_{{\rm{LTI}}} (\hat {\cal
F}_{{\rm{SF}}}^{(N,\, \downarrow )} ,\,\Sigma ) - P\nonumber\\
&\quad = {\cal H}(\hat {\cal F}_{{\rm{SF}}}^{(N,\, \downarrow )}
)_ \bot ^T \Pi _N^T {\cal X}_N (P,\,1)\Pi _N {\cal H}(\hat {\cal
F}_{{\rm{SF}}}^{(N,\, \downarrow )} )_ \bot\nonumber\\
&\quad \prec 0,\label{theorem2-eq2}
\end{align}
where
\begin{align*}
&{\cal H}(\hat {\cal F}_{{\rm{SF}}}^{(N,\, \downarrow )} )_ \bot =
\left[ {\begin{array}{*{20}c}
   {({\bf{I}}_N  \otimes {\bf{I}}_n )}  \\
   {\hat {\cal F}_{{\rm{SF}}}^{(N,\, \downarrow )} }  \\
\end{array}} \right]\left[ {\begin{array}{*{20}c}
   {{\bf{I}}_n }  \\
   {{\cal A}_{{\rm{LTI}}} (\hat {\cal F}_{{\rm{SF}}}^{(1,\, \downarrow )} ,\,\Sigma )}  \\
    \vdots   \\
   {{\cal A}_{{\rm{LTI}}} (\hat {\cal F}_{{\rm{SF}}}^{(N - 1,\, \downarrow )} ,\,\Sigma )}  \\
\end{array}} \right].
\end{align*}
Then, relying on the elimination lemma \citep{Boyd1994}, we prove
that both \eqref{theorem2-eq1} and \eqref{theorem2-eq2} are
satisfied if and only if there exists ${\cal V}$ such that
\eqref{theorem1-eq2} holds with $M = {\cal H}(\hat {\cal
F}_{{\rm{SF}}}^{(N,\, \downarrow )} )^T {\cal V}$. This proves the
sufficiency.

(Necessity) Assume that BMIs \eqref{theorem1-eq1} and
\eqref{theorem1-eq2} with $M = {\cal H}(\hat {\cal
F}_{{\rm{SF}}}^{(N,\, \downarrow )} )^T {\cal V}$ admit a solution
$(P,\,{\cal F}_{{\rm{SOF}}}^{(N,\, \downarrow )} )$. By means of
the elimination lemma, we have that \eqref{theorem2-eq1} and
\eqref{theorem2-eq2} hold. This implies $P \in {\cal P}({\cal
A}_{{\rm{LTI}}} (\hat {\cal F}_{{\rm{SF}}}^{(N,\, \downarrow )}
,\,\Sigma ))$ and $P \in {\cal S}_{{\rm{SOF}}}^{(N)} (\Sigma )$,
so \eqref{nonempty-condition} is satisfied. This completes the
proof. $\blacksquare$

Theorem \ref{theorem2} tells us that if $(\hat P,\,\hat {\cal J}
,\,\hat {\cal G})$ is a solution to \eqref{Ebihara2011-lemma-eq1},
then $M = {\cal H}(\hat {\cal F}_{{\rm{SF}}}^{(N,\, \downarrow )}
)^T {\cal V}$ with appropriately selected ${\cal V}$ can be a
reasonable choice of $M$ so that \eqref{theorem1-eq1} and
\eqref{theorem1-eq2} in Theorem \ref{theorem1} become feasible.
The following corollary can be immediately obtained from Theorems
\ref{theorem1} and \ref{theorem2}:
\begin{corollary}\label{corollary1}
Suppose that
\begin{enumerate}
\item[a)] $(\hat P,\,\hat {\cal J} ,\,\hat {\cal G})$ is a
solution to \eqref{Ebihara2011-lemma-eq1}, and $\hat {\cal
F}_{{\rm{SF}}}^{(N,\, \downarrow )}  = ({\bf{T}}_N \otimes
{\bf{I}}_m )\hat {\cal J} \hat {\cal G}^{-1} ({\bf{T}}_N \otimes
{\bf{I}}_n )$;

\item[b)] \eqref{nonempty-condition} is satisfied.
\end{enumerate}
Then, there exists ${\cal F}_{{\rm{SOF}}}^{(N,\, \downarrow )} \in
{\mathbb R}^{Nm \times Np}$ defined in \eqref{Fsof} such that
$N$-periodic control system \eqref{SOF-closed-loop-system} or
equivalent LTI system \eqref{SOF-LTI-system} is asymptotically
stable if and only if there exist matrices $P = P^T \in {\mathbb
R}^{n \times n} $ and ${\cal V} \in {\mathbb R}^{((N - 1)n + Nm)
\times ((N - 1)n + Nm)}$ such that the following problem is
satisfied with ${\cal F}_{{\rm{SOF}}}^{(N,\, \downarrow )}  \in
{\mathbb R}^{Nm \times Np}$:
\begin{align}
&P \succ 0,\label{corollary1-eq1}\\
& \Pi _N^T {\cal X}_N (P,\,1)\Pi _N  + {\rm{He}}\{ {\cal H}(\hat
{\cal F}_{{\rm{SF}}}^{(N,\, \downarrow )} )^T {\cal V}{\cal
C}({\cal F}_{{\rm{SOF}}}^{(N,\, \downarrow )} )\}  \prec
0.\label{corollary1-eq2}
\end{align}
\end{corollary}
{\bf Proof.} The sufficiency follows immediately from Theorem
\ref{theorem1}. To prove the necessity, suppose that there exists
${\cal F}_{{\rm{SOF}}}^{(N,\, \downarrow )}$ defined in
\eqref{Fsof} such that $N$-periodic control system
\eqref{SOF-closed-loop-system} is asymptotically stable. Since
\eqref{nonempty-condition} is satisfied by assumption, one can
select ${\cal F}_{{\rm{SOF}}}^{(N,\, \downarrow )}$ so that ${\cal
P}({\cal A}_{{\rm{LTI}}} (\hat {\cal F}_{{\rm{SF}}}^{(N,\,
\downarrow )} ,\,\Sigma )) \cap {\cal P}({\cal A}_{{\rm{LTI}}}
({\cal F}_{{\rm{SOF}}}^{(N,\, \downarrow )} ,\,\Sigma )) \ne
\emptyset$. The rest of the proof then follows the same line as in
the sufficient part of the proof of Theorem \ref{theorem2}.
$\blacksquare$
\begin{remark}
In the case $N=1$, \eqref{corollary1-eq1} and
\eqref{corollary1-eq2} reduce to Theorem 1 in \cite{Peaucelle2001}
and \cite{Arzelier2010}.
\end{remark}

It should be kept in mind that there is no guarantee that $\hat
{\cal F}_{{\rm{SF}}}^{(N,\, \downarrow )} \in {\cal
L}_{{\rm{SF}}}^{(N)}$ obtained by solving LMI
\eqref{Ebihara2011-lemma-eq1} satisfies condition
\eqref{nonempty-condition}. In addition, since
\eqref{nonempty-condition} is used only in the necessity part of
the proof of Corollary \ref{corollary1}, in a practical
implementation, Corollary \ref{corollary1} should be regarded as
only a sufficient condition. Note also that the condition of
Corollary \ref{corollary1} is still a BMI problem. However, as in
\cite{Peaucelle2001}, one can expect that solving the BMI of
Corollary \ref{corollary1} gives better results than solving the
BMI of Theorem \ref{theorem1}, since an initial guess of $M$ in
Theorem \ref{theorem1} is used in Corollary \ref{corollary1}.
Unfortunately, if $\hat {\cal F}_{{\rm{SF}}}^{(N,\, \downarrow )}
\in {\cal L}_{{\rm{SF}}}^{(N)}$ does not satisfy
\eqref{nonempty-condition}, the BMI of Corollary \ref{corollary1}
has no solution even when the solution set of the original BMI of
Theorem \ref{theorem1} is nonempty. In this respect, it can be
said that another source of conservatism is introduced in the BMI
of Corollary \ref{corollary1}. Conceptually, we conjecture that
this conservatism can be reduced by increasing $N$. To give an
intuitive perspective on how increasing $N$ reduces this kind of
conservatism, let us introduce the following lemmas:
\begin{lemma}\label{basic-lamma1}
Assume ${\mathcal S}_{{\rm{SF}}}^{(1)} (\Sigma ) \ne \emptyset $.
The following statements are true:
\begin{enumerate}[a)]
\item ${\mathcal S}_{{\rm{SF}}}^{(1)} (\Sigma ) \subseteq
{\mathcal S}_{{\rm{SF}}}^{(N)} (\Sigma ),\quad \forall N \in
{\mathbb N}_+$.

\item $\mathop {\lim }\limits_{N \to \infty } {\cal
S}_{{\rm{SF}}}^{(N)} (\Sigma ) = {\mathbb S}_+^n$.
\end{enumerate}
\end{lemma}
{\bf Proof.} a) For any $P \in {\mathcal S}_{{\rm{SF}}}^{(1)}
(\Sigma )$, assume that $F = {\mathcal F}_{{\rm{SF}}}^{(1,\,
\downarrow )}$ satisfies ${\mathcal A}_{{\rm{LTI}}} ({\mathcal
F}_{{\rm{SF}}}^{(1,\, \downarrow )} ,\,\Sigma )^T P{\mathcal
A}_{{\rm{LTI}}} ({\mathcal F}_{{\rm{SF}}}^{(1,\, \downarrow )}
,\,\Sigma ) - P = (A + BF)^T P(A + BF) - P \prec 0$. Then, $P \in
{\mathcal S}_{{\rm{SF}}}^{(N)} (\Sigma )$ because $(A + BF)^{NT}
P(A + BF)^N  - P = {\cal A}_{{\rm{LTI}}} ({\cal
F}_{{\rm{SF}}}^{(N,\, \downarrow )} ,\,\Sigma )^T P{\cal
A}_{{\rm{LTI}}} ({\cal F}_{{\rm{SF}}}^{(N,\, \downarrow )}
,\,\Sigma ) - P \prec 0$ holds with ${\mathcal
F}_{{\rm{SF}}}^{(N,\, \downarrow )}  = {\bf{I}}_N \otimes F$. This
implies a) is true.

b) For any ${\mathcal F}_{{\rm{SF}}}^{(1,\, \downarrow )}  \in
{\mathcal L}_{{\rm{SF}}}^{(1)}$ and $P \in {\mathbb S}_+^n$, it
holds that $\lim _{N \to \infty } ({\mathcal A}_{{\rm{LTI}}}
({\mathcal F}_{{\rm{SF}}}^{(1,\, \downarrow )} ,\,\Sigma )^{NT}
P{\mathcal A}_{{\rm{LTI}}} ({\mathcal F}_{{\rm{SF}}}^{(1,\,
\downarrow )} ,\,\Sigma )^N  - P) =  - P \prec 0$. Since
${\mathcal A}_{{\rm{LTI}}} ({\mathcal F}_{{\rm{SF}}}^{(1,\,
\downarrow )} ,\,\Sigma )^N  = {\mathcal A}_{{\rm{LTI}}}
({\mathcal F}_{{\rm{SF}}}^{(N,\, \downarrow )} ,\,\Sigma )$ with
${\mathcal F}_{{\rm{SF}}}^{(N,\, \downarrow )} = {\bf{I}}_N
\otimes {\mathcal F}_{{\rm{SF}}}^{(1,\, \downarrow )}$ , we have
\[\lim _{N \to \infty } ({\cal A}_{{\rm{LTI}}} ({\cal
F}_{{\rm{SF}}}^{(N,\, \downarrow )} ,\,\Sigma )^T P{\cal
A}_{{\rm{LTI}}} ({\cal F}_{{\rm{SF}}}^{(N,\, \downarrow )}
,\,\Sigma ) - P) \prec 0\] for any $P \in {\mathbb S}_+^n $. This
implies b) is satisfied, and the proof is completed.
$\blacksquare$
\begin{lemma}\label{basic-lamma2}
Assume ${\mathcal S}_{{\rm{SF}}}^{(1)} (\Sigma ) \ne \emptyset $.
The following statements are true:
\begin{enumerate}[a)]
\item ${\cal S}_{{\rm{SOF}}}^{(1)} (\Sigma ) \subseteq {\cal
S}_{{\rm{SOF}}}^{(N)} (\Sigma ),\quad \forall N \in {\mathbb N}_ +
$.

\item $\mathop {\lim }\limits_{N \to \infty } {\cal
S}_{{\rm{SOF}}}^{(N)} (\Sigma ) = {\mathbb S}_+^n $.

\item ${\cal S}_{{\rm{SOF}}}^{(N)} (\Sigma ) \subseteq {\cal
S}_{{\rm{SF}}}^{(N)} (\Sigma ),\quad \forall N \in {\mathbb N}_+$

\end{enumerate}
\end{lemma}
{\bf Proof.} Proofs for statements a) and b) follow immediately
from those of Lemma 2. For statement c), assume that $P \in {\cal
S}_{{\rm{SOF}}}^{(N)} (\Sigma )$, which means there exists ${\cal
F}_{{\rm{SOF}}}^{(N,\, \downarrow )}  \in {\cal
L}_{{\rm{SOF}}}^{(N)}$ satisfying ${\cal A}_{{\rm{LTI}}} ({\cal
F}_{{\rm{SOF}}}^{(N,\, \downarrow )} ,\,\Sigma )^T P{\cal
A}_{{\rm{LTI}}} ({\cal F}_{{\rm{SOF}}}^{(N,\, \downarrow )}
,\,\Sigma ) - P \prec 0$. Then, $P \in {\cal S}_{{\rm{SF}}}^{(N)}
(\Sigma )$ because
\begin{align*}
&{\cal A}_{{\rm{LTI}}} ({\cal F}_{{\rm{SF}}}^{(N,\, \downarrow )}
,\,\Sigma )^T P{\cal A}_{{\rm{LTI}}} ({\cal F}_{{\rm{SF}}}^{(N,\,
\downarrow )} ,\,\Sigma ) - P \prec 0
\end{align*}
holds with ${\cal F}_{{\rm{SF}}}^{(N,\, \downarrow )}  = {\cal
F}_{{\rm{SOF}}}^{(N,\, \downarrow )} ({\bf{I}}_N  \otimes C)$.
This completes the proof. $\blacksquare$

Again, recall that $(\hat P,\,\hat {\cal J},\,\hat {\cal G})$ is a
solution to \eqref{Ebihara2011-lemma-eq1}, and $\hat {\cal
F}_{{\rm{SF}}}^{(N,\, \downarrow )}  = ({\bf{T}}_N \otimes
{\bf{I}}_m )\hat {\cal J} \hat {\cal G}^{- 1} ({\bf{T}}_N \otimes
{\bf{I}}_n ) \in {\cal L}_{{\rm{SF}}}^{(N)}$. Let us assume ${\cal
S}_{{\rm{SOF}}}^{(1)} (\Sigma ) \ne \emptyset$. Then, in view of
Lemmas \ref{basic-lamma1} and \ref{basic-lamma2}, it is true that
${\cal S}_{{\rm{SOF}}}^{(N)} (\Sigma ) \subseteq {\cal
S}_{{\rm{SF}}}^{(N)} (\Sigma ),\, \forall N \in {\mathbb N}_+$,
$\lim _{N \to \infty } {\cal S}_{{\rm{SF}}}^{(N)} (\Sigma ) \cap
{\cal S}_{{\rm{SOF}}}^{(N)} (\Sigma ) = \lim _{N \to \infty }
{\cal S}_{{\rm{SOF}}}^{(N)} (\Sigma ) = {\mathbb S}_+^n$, and
$\lim _{N \to \infty } \{ P \in {\mathbb S}_ + ^n :\,P \in {\cal
S}_{{\rm{SF}}}^{(N)} (\Sigma ),\,P \notin {\cal
S}_{{\rm{SOF}}}^{(N)} (\Sigma )\}  = \emptyset$. In addition, let
us suppose that $\hat P^{ - 1}  \in {\cal P}({\cal A}_{{\rm{LTI}}}
(\hat {\cal F}_{{\rm{SF}}}^{(N,\, \downarrow )} ,\,\Sigma ))$ is a
random matrix within ${\cal S}_{{\rm{SF}}}^{(N)} (\Sigma )$. Then,
we can expect that as $N$ gets larger, set $\{ P \in {\mathbb S}_
+ ^n :\,P \in {\cal S}_{{\rm{SF}}}^{(N)} (\Sigma ),\,P \notin
{\cal S}_{{\rm{SOF}}}^{(N)} (\Sigma )\}$ tends to shrink and
eventually become the empty set as $N \to \infty $. Thus, it is
more likely that $\hat P^{ - 1} \in {\cal P}({\cal A}_{{\rm{LTI}}}
(\hat {\cal F}_{{\rm{SF}}}^{(N,\, \downarrow )} ,\,\Sigma ))
\subseteq {\cal S}_{{\rm{SF}}}^{(N)} (\Sigma )$ lies within ${\cal
S}_{{\rm{SOF}}}^{(N)} (\Sigma ) \cap {\cal S}_{{\rm{SF}}}^{(N)}
(\Sigma )$ as $N \to \infty $. In other words, as $N$ increases,
there is a more possibility that \eqref{nonempty-condition} holds,
and thus, the solution set of the BMI problem of Corollary
\ref{corollary1} is nonempty.

To determine ${\cal F}_{{\rm{SOF}}}^{(N,\, \downarrow )} $, the
problem of Corollary \ref{corollary1} is still a BMI problem (not
an LMI in ${\cal V}$ and ${\cal F}_{{\rm{SOF}}}^{(N,\, \downarrow
)}$). Local solutions to the BMI problem of Corollary
\ref{corollary1} can be obtained by using PENBMI
\citep{Kocvara2005}. Alternatively, with a suitable choice of
particular ${\cal V}$, the BMI can reduce to a convex LMI at the
price of some conservatism. For instance, letting
\begin{align}
&{\cal V} = \left[ {\begin{array}{*{20}c}
   {{\cal V}_{11}} & {{\cal V}_{12}}  \\
   {\bf{0}} & {{\cal V}_{22}}  \\
\end{array}} \right],\label{V-matrix}
\end{align}
where
\begin{align}
&{\cal V}_{11} : = \left[ {\begin{array}{*{20}c}
   {V_{11}^{(1,\,1)} } & {\bf{0}} & {\bf{0}}  \\
    \vdots  &  \ddots  &  \ddots   \\
   {V_{11}^{(N,\,1)} } &  \cdots  & {V_{11}^{(N,\,N)} }  \\
\end{array}} \right] \in {\mathbb R}^{Nm \times Nm},\label{V11}
\end{align}
$V_{11}^{(i,\,j)}  \in {\mathbb R}^{m \times m} ,\,{\cal V}_{12}
\in {\mathbb R}^{Nm \times (N - 1)n}$, and ${\cal V}_{22}\in
{\mathbb R}^{(N - 1)n \times (N - 1)n}$, the following result is
obtained:
\begin{corollary}\label{corollary3}
Suppose that $(\hat P,\,\hat {\cal J},\,\hat {\cal G})$ is a
solution to \eqref{Ebihara2011-lemma-eq1}, and let $\hat {\cal
F}_{{\rm{SF}}}^{(N,\, \downarrow )}  = ({\bf{T}}_N \otimes
{\bf{I}}_m )\hat {\cal J}\hat {\cal G}^{- 1} ({\bf{T}}_N \otimes
{\bf{I}}_n )$. Then, system \eqref{system} is stabilizable via
$N$-PTVMSOFC \eqref{PTVMSOFC} if there exists matrices $P = P^T
\in {\mathbb R}^{n \times n} ,\,M^{(i,\,j)} \in {\mathbb R}^{m
\times p} ,\,V_{11}^{(i,\,j)}  \in {\mathbb R}^{m \times
m},\,{\cal V}_{12}  \in {\mathbb R}^{Nm \times (N - 1)n}$, and
${\cal V}_{22}\in {\mathbb R}^{(N - 1)n \times (N - 1)n}$ such
that the following LMI problem is satisfied with $\gamma  = 1$:
\begin{align}
&P \succ 0,\label{corollary3-eq1}\\
&\Pi _N^T {\cal X}_{N} (P,\,\gamma )\Pi _N \nonumber\\
&\quad + {\rm{He}}\{ {\cal H}(\hat {\cal F}_{{\rm{SF}}}^{(N,\,
\downarrow )} )^T {\cal D}({\cal M},\,{\cal V}_{11} ,\,{\cal
V}_{12} ,\,{\cal V}_{22} )\}
\prec 0,\label{corollary3-eq2}\\
&{\rm{He}}\{ {\cal V}_{11} \}  \prec 0,\label{corollary3-eq3}
\end{align}
where ${\cal V}_{11}$ is defined in \eqref{V11},
\begin{align}
&{\cal M}: = \left[ {\begin{array}{*{20}c}
   {M^{(1,\,1)} } & {\bf{0}} & {\bf{0}}  \\
    \vdots  &  \ddots  &  \ddots   \\
   {M^{(N,\,1)} } &  \cdots  & {M^{(N,\,N)} }  \\
\end{array}} \right] \in {\mathbb R}^{Nm \times Np},\label{M}
\end{align}
and ${\cal D}({\cal M},\,{\cal V}_{11} ,\,{\cal V}_{12} ,\,{\cal
V}_{22} )$ is defined in \eqref{corollary3-eq4} at the top of the
next page. Moreover, an admissible $N$-PTVMSOFC gain matrix is
given by ${\cal F}_{{\rm{SOF}}}^{(N,\, \downarrow )}  = {\cal
V}_{11}^{- 1} {\cal M}$.
\begin{figure*}[!t]
\begin{align}
&{\cal D}({\cal M},\,{\cal V}_{11} ,\,{\cal V}_{12} ,\,{\cal
V}_{22} ): = \left[ {\begin{array}{*{20}c}
   {{\cal M}({\bf{I}}_N  \otimes C) + {\cal V}_{12} ({\cal L}_{N - 1}  \otimes A) - {\cal V}_{12} ({\cal R}_{N - 1}  \otimes {\bf{I}}_n )} &\vline &  { - {\cal V}_{11}  + {\cal V}_{12} ({\cal L}_{N - 1}  \otimes B)}  \\
\hline {{\cal V}_{22} ({\cal L}_{N - 1}  \otimes A) - {\cal V}_{22} ({\cal R}_{N - 1}  \otimes {\bf{I}}_n )} &\vline &  {{\cal V}_{22} ({\cal L}_{N - 1}  \otimes B)}  \\
\end{array}} \right].\label{corollary3-eq4}
\end{align}
\hrulefill \vspace*{4pt}
\end{figure*}
\end{corollary}
{\bf Proof.} Noting that \eqref{corollary3-eq3} ensures the
invertibility of ${\cal V}_{11}$, substituting \eqref{V-matrix}
into \eqref{corollary1-eq2}, and using the change of variables
${\cal M} = {\cal V}_{11} {\cal F}_{{\rm{SOF}}}^{(N,\, \downarrow
)}$, we have that \eqref{corollary3-eq2} is equivalent to
\eqref{corollary1-eq2}. $\blacksquare$
\begin{remark}
LMI \eqref{corollary3-eq3} guarantees that ${\cal V}_{11}$ is
nonsingular. If it is eliminated, then LMIs \eqref{corollary3-eq1}
and \eqref{corollary3-eq2} can yield less conservative results
although the invertibility of ${\cal V}_{11}$ is not guaranteed.
Therefore, instead of
\eqref{corollary3-eq1}-\eqref{corollary3-eq3}, we can use only
\eqref{corollary3-eq1} and \eqref{corollary3-eq2}, and when they
are feasible, the invertibility of ${\cal V}_{11}$ should be
checked to ultimately determine the feasibility of the control
design problem.
\end{remark}

Based on Corollary \ref{corollary3}, the two-steps algorithms
suggested in \cite{Mehdi2004,Agulhari2010, Agulhari2012} can be
adopted to design the $N$-PTVMSOFC.

\vspace*{-12pt}\hrulefill\vspace*{-12pt}

{\bf Algorithm 2. Two-Steps LMI Algorithm.}

\vspace*{-17pt}\hrulefill\vspace*{-5pt}
\begin{enumerate}[Step 1.]

\item Solve LMI \eqref{Ebihara2011-lemma-eq1} for $(P,\,{\cal
G},\,{\cal J})$ and let $\hat {\cal F}_{{\rm{SF}}}^{(N,\,
\downarrow )}  = ({\bf{T}}_N  \otimes {\bf{I}}_m )\hat {\cal J}
\hat {\cal G}^{- 1} ({\bf{T}}_N  \otimes {\bf{I}}_n )$ with $(\hat
P,\,\hat {\cal G},\,\hat {\cal J}) \in \{ P,\,{\cal G} ,\,{\cal
J}:\,{\rm{LMI}}\,\eqref{Ebihara2011-lemma-eq1}\}$.

\item With $\hat {\cal F}_{{\rm{SF}}}^{(N,\, \downarrow )}$
obtained from the previous step, solve for  $\Lambda : =
(P,\,{\cal M},\,{\cal V}_{11},\,{\cal V}_{12},\,{\cal V}_{12})$
LMIs \eqref{corollary3-eq1}-\eqref{corollary3-eq3} with
$\gamma=1$:
\begin{align*}
&(\hat P,\,\hat {\cal M},\,\hat {\cal V}_{11},\,\hat
{\cal V}_{12},\,\hat {\cal V}_{22})\\
&\quad \in \{ \Lambda
:\,{\rm{LMIs}}\,\eqref{corollary3-eq1}-\eqref{corollary3-eq3},\,
\gamma=1\}.
\end{align*}

If feasible, then $\hat {\cal F}_{{\rm{SOF}}}^{(N,\, \downarrow )}
= \hat {\cal V}_{11}^{- 1} \hat {\cal M}$ is a stabilizing
$N$-PTVMSOFC gain matrix.
\end{enumerate}
\vspace*{-17pt}\hrulefill \vspace*{0pt}

Moreover, versions of the ILMI algorithm that alternates between
the SF and the SOF designs developed in \cite{Peaucelle2001} can
be also applied as less conservative alternatives.

\vspace*{-12pt}\hrulefill\vspace*{-12pt}

 {\bf Algorithm 3. ILMI Algorithm.}

\vspace*{-17pt}\hrulefill\vspace*{-5pt}
\begin{enumerate}[Step 1.]
\item (Initialization). Set $i=1$, the maximum number of
iterations $N_{{\rm{iter}}}  \in {\mathbb N}_+$, and a
sufficiently small positive real number $\delta$. Solve LMI
\eqref{Ebihara2011-lemma-eq1} for $(P,\,{\cal G},\,{\cal J})$ and
let $\hat {\cal F}_{{\rm{SF}}}^{(N,\, \downarrow )}  = ({\bf{T}}_N
\otimes {\bf{I}}_m )\hat {\cal J}\hat {\cal G}^{- 1} ({\bf{T}}_N
\otimes {\bf{I}}_n )$ with $(\hat P,\,\hat {\cal G},\,\hat {\cal
J}) \in \{ P,\,{\cal G},\,{\cal
J}:\,{\rm{LMI}}\,\eqref{Ebihara2011-lemma-eq1}\}$.

\item With $\hat {\cal F}_{{\rm{SF}}}^{(N,\, \downarrow )}$
obtained from the previous step, solve for $\Lambda : = (\gamma
,\,P,\,{\cal M},\,{\cal V}_{11} ,\,{\cal V}_{12},\,{\cal V}_{22})$
the following optimization problem:
\begin{align}
&(\hat \gamma_i ,\,\hat P,\,\hat {\cal M},\,\hat {\cal V}_{11}
,\,\hat {\cal V}_{12},\,\hat {\cal
V}_{22})\nonumber\\
&\quad : = \arg \mathop {\min }\limits_\Lambda  \{ \gamma  \in
{\mathbb
R}:\,{\rm{LMIs}}\,\eqref{corollary3-eq1}-\eqref{corollary3-eq3}\}.\label{optimization1}
\end{align}

\item If $\hat \gamma _i  \le 1$, then $\hat {\cal
F}_{{\rm{SOF}}}^{(N,\, \downarrow )}  = \hat {\cal V}_{11}^{- 1}
\hat {\cal M}$ is a stabilizing $N$-PTVMSOFC gain matrix. STOP.
Otherwise, if $i \ge 2$ and $|\hat \gamma _{i - 1}  - \hat \gamma
_i | \le \delta$ or $\hat \gamma _{i - 1} < \hat \gamma _i$ or $i
= N_{{\rm{iter}}} $, then this algorithm cannot get a feasible
solution. STOP.

\item With $(\hat \gamma _i ,\,\hat {\cal M},\,\hat {\cal
V}_{11},\,\hat {\cal V}_{12},\,\hat {\cal V}_{22})$ obtained from
the previous step, solve the LMI problem
\begin{align*}
&(\hat P,\,\hat {\cal F}_{{\rm{SF}}}^{(N,\, \downarrow )} )\\
&\in \left\{ \begin{array}{l}
 (P,\,{\cal F}_{{\rm{SF}}}^{(N,\, \downarrow )} ):\,P \succ 0,\,\Pi _N^T {\cal X}_N (P,\,\hat \gamma _i )\Pi _N  \\
  + {\rm{He}}\{ {\cal H}({\cal F}_{{\rm{SF}}}^{(N,\, \downarrow )} )^T {\cal D}(\hat {\cal M},\,\hat {\cal V}_{11} ,\,\hat {\cal V}_{12} ,\,\hat {\cal V}_{22} )\}  \\
  \prec 0 \\
 \end{array} \right\},
\end{align*}
set $i = i + 1$, and go to Step 2.
\end{enumerate}
\vspace*{-17pt}\hrulefill \vspace*{0pt}
\begin{remark}
The optimization problem \eqref{optimization1} is a unidimensional
minimization subject to LMI constraints, and for fixed $\gamma$,
conditions \eqref{corollary1-eq1} and \eqref{corollary1-eq2} are
LMIs tractable via LMI solvers
\citep{Gahinet1995,Lofberg2004,Strum1999}. Thus, the optimization
problem can be solved by means of a sequence of LMI problems, i.e.
a line search or a bisection process over $\gamma$. Moreover, the
optimization problem belongs to the class of eigenvalue problems,
which are convex optimizations \citep{Boyd1994}, and hence, can be
directly treated with the aid of the LMI solver
\citep{Gahinet1995}.
\end{remark}
\begin{remark}
It is not difficult to show that, at least theoretically, if the
LMI problem at Step 2 and $i = 1$ is feasible, then all the
subsequent LMIs are also feasible for all $i > 1$, and $\{ \hat
\gamma _1 ,\,\hat \gamma _2 , \ldots \} $ is a conversing and
non-increasing sequence. However, in practice, the LMIs after Step
2 can fail to find a feasible solution or $\hat \gamma _i$ can
increase and fluctuate irregularly in many cases. This phenomenon
may be common to many other ILMI schemes and may be due to the
fact that as solution spaces of the LMIs become narrower, the
feasibility of the LMIs tends to be more sensitive to small
numerical errors of the solutions computed at the previous steps.
In this case, the algorithm can be deemed not to be able to get a
solution.
\end{remark}

All computations in the sequel were done in MATLAB R2012b running
under Windows 7 PC. The computer used was equipped with an Intel
Core i7-3770 3.4GHz CPU and 32GB RAM. The LMI problems were solved
with SeDuMi \citep{Strum1999} and Yalmip \citep{Lofberg2004}.

\begin{example}\label{example1}
For a statistical comparison analysis of the proposed results with
existing ones, we randomly generated thousand systems with
$(n,\,m,\,p) = (3,\,1,\,1)$ whose open-loop systems were unstable.
Each system was computed using the following procedure: 1) triplet
$(A,\,B,\,C)$ is generated with matrices whose entries are real
numbers uniformly distributed in the interval $[ - 2,\,2]$; 2) $A$
is replaced with $(1.2/\rho (A))A$ so that the spectral radius of
$A$ becomes $1.2$; 3) if $(A,\,B)$ is stabilizable and $(C,\,A)$
is detectable, then add the triplet to the list of test systems.
Else, discard it and go to step 1). Since the PTVMSOFC can be
interpreted as a sort of dynamic output feedback (DOF) controller,
the proposed approaches are also compared with the full-order DOF
design \citep{Iwasaki1994,Scherer1997}. The number of stabilizable
systems, denoted by $N_{{\rm{stable}}}$, in the context of
feasibility of several approaches are listed in Table \ref{table1}
with the average computational time (in seconds) spent by each
test, where for Algorithm 3, we set $(N_{iter} ,\,\delta ) =
(10,\,10^{ - 4} )$, and for optimization \eqref{optimization1}, a
bisection algorithm over $\gamma $ was used. In addition, for
PENBMI, we used the BMI condition
\begin{align*}
&\left[ {\begin{array}{*{20}c}
   { - P} & {(A + BFC)^T P}  \\
   {P(A + BFC)} & { - P}  \\
\end{array}} \right] \prec 0,
\end{align*}
From Table \ref{table1}, the following observation can be made:
\begin{enumerate}[a)]
\item The results show that at the price of a higher computational
cost, the proposed method offers improvement over the previous
approaches except for the full-order DOF design. The number of
parameters of the controller is $mpN(N + 1)/2$ for the PTVMSOFC
while $n^2  + np + mn + mp$ for the full-order DOF controller.

In order to compare and evaluate the on-line computational burden,
we will check the number of operations including multiplication
and addition. The total multiplication and addition during period
$k \in \{ 0,\,1,\, \ldots ,\,N - 1\} $ are summarized in Table
\ref{FODOF-comp-cost} for the full-order DOF and Table
\ref{PTVMSOFC-comp-cost} for the PTVMSOFC.
\begin{table}[t]
\caption{Full-order DOF, period $k \in \{ 0,\,1,\, \ldots ,\,N -
1\}$}
\begin{center}
\begin{tabular}{ c c}
\hline Multiplication & Addition \\
\hline $N(n^2  + np + mn + mp)$   & $N(m + n)(n + p - 1)$\\
\hline
\end{tabular}
\end{center}
\label{FODOF-comp-cost}
\end{table}

\begin{table}[t]
\caption{PTVMSOFC, period $k \in \{ 0,\,1,\, \ldots ,\,N - 1\}$}
\begin{center}
\begin{tabular}{ c c}
\hline Multiplication & Addition \\
\hline $mpN(N + 1)/2$ & $mpN(N + 1)/2 - Nm$
\\
\hline
\end{tabular}
\end{center}
\label{PTVMSOFC-comp-cost}
\end{table}

It might not be an easy task to perform the qualitative analysis
for a large number of combinations of $(N,\,p,\,m,\,n)$. However,
by investigating a simple example, we can observe that the
off-line computation of the PTVMSOFC can be smaller than that of
the full-order DOF in some cases. The examples are shown in Table
\ref{FODOF-comp-cost-ex} for the full-order DOF and Table
\ref{PTVMSOFC-comp-cost-ex} for the PTVMSOFC.
\begin{table}[t]
\caption{Full-order DOF, period $k \in \{ 0,\,1,\, \ldots ,\,N -
1\}$}
\begin{center}
\begin{tabular}{c c c}
\hline $n$ &  Multiplication & Addition \\
  & $N(n^2  + np + mn + mp)$   & $N(m + n)(n + p - 1)$\\
\hline  2 & 18 & 12 \\
\hline  3 & 32 & 24 \\
\hline  4 & 50 & 40 \\
\hline  5 & 72 & 60 \\
\hline  6 & 98 & 84 \\
\hline
\end{tabular}
\end{center}
\label{FODOF-comp-cost-ex}
\end{table}
\begin{table}[t]
\caption{PTVMSOFC, period $k \in \{ 0,\,1,\, \ldots ,\,N - 1\}$}
\begin{center}
\begin{tabular}{c c c}
\hline $n$ & Multiplication & Addition \\
 & $mpN(N + 1)/2$ & $mpN(N + 1)/2 - Nm$\\

\hline 2 & 3 & 1 \\
\hline 3 & 3 & 1 \\
\hline 4 & 3 & 1 \\
\hline 5 & 3 & 1 \\
\hline 6 & 3 & 1 \\
\hline
\end{tabular}
\end{center}
\label{PTVMSOFC-comp-cost-ex}
\end{table}

By comparing the results, it can be seen that the on-line
computational cost of the PTVMSOFC can be lower than that of the
full-order DOF in some cases. In the above case, it is interesting
to observe that the computational cost of the PTVMSOFC is not
dependent on dimension $n$ of the state.

In summary, the online computational cost of the PTVMSOFC can be
lower than that of the full-order DOF. For this reason, the
PTVMSOFC can be a useful alternative to the DOF controller in some
cases.

\item It can be observed that solving Corollary \ref{corollary1}
with PENBMI is less conservative than solving Theorem
\ref{theorem1} with PENBMI. Since Corollary \ref{corollary1} is
derived from Theorem \ref{theorem1} with a reasonable initial
selection of $M$ based on the SF design, we can conclude that the
improvement of Corollary \ref{corollary1} mainly comes from the
initialization of $M$. Moreover, the comparison results between
the two-steps algorithm and the ILMI algorithm suggest that some
improvement can be achieved by adopting the ILMI method.
\end{enumerate}

\begin{table*}[t]
\caption{Example \ref{example1}. Number of stabilizable systems,
$N_{{\rm{stable}}}$, and the average computational time.}
\begin{center}
\begin{tabularx}{\textwidth}{ p{14cm} Z Z}
\hline

\centering  Methods & $N_{{\rm{stable}}}$ & Time (s)\\

\hline

\rowcolor[gray]{.9}\centering Cone complementarity linearization algorithm in \cite{Ghaoui1997} (discrete-time version) & $509$ & $34.35$\\
\rowcolor[gray]{.9}\centering Discrete $P$-problem in
\cite{Crusius1999} & $248$ & $0.10$\\
\rowcolor[gray]{.9}\centering Discrete $W$-problem in \cite{Crusius1999} & $243$ & $0.10$\\

\rowcolor[gray]{.9}\centering Algorithm A in \cite{Rosinova2003}
with $(R,\,Q) = (0.01{\bf{I}}_m ,\,{\bf{I}}_n )$ & $306$ & $0.04$\\

\rowcolor[gray]{.9}\centering Two-steps LMI approach of Theorem 3.1 in \cite{Mehdi2004} with constraint $-G-G^T \prec 0$ & $427$ & $0.18$\\

\rowcolor[gray]{.9}\centering Two-steps LMI approach of Theorem 3.1 in \cite{Mehdi2004} with $F_1=F_3= {\bf 0}$ & $427$ & $0.18$\\

\rowcolor[gray]{.9}\centering Lemma 3 in \cite{Dong2007} with $T = [C^T (CC^T )^{ - 1} \,\,\,C_ \bot  ]$ (Method in \cite{de_Oliveira2002}) & $287$ & $0.10$\\

\rowcolor[gray]{.9}\centering Theorems 3.1 and 3.3 in
\cite{Bara2005} with $T = [C^T (CC^T )^{ - 1} \,\,\,C_ \bot  ]$ & $222$ & $0.14$\\

\rowcolor[gray]{.9}\centering Algorithm 1 in \cite{Shu2010} & $309 $ & $9.90$\\

\rowcolor[gray]{.9}\centering PENBMI \citep{Kocvara2005} & $484$ & $0.08$ \\

\rowcolor[gray]{.9}\centering Full-order DOF design (discrete-time version of \cite{Scherer1997})  & $1000$ & $0.10$ \\

\rowcolor[gray]{.85}\centering Theorem \ref{theorem1} solved with PENBMI for $N=1$ & $355$ & $0.09$\\
\rowcolor[gray]{.85}\centering Theorem \ref{theorem1} solved with PENBMI for $N=2$ & $495$ & $0.25$\\
\rowcolor[gray]{.85}\centering Theorem \ref{theorem1} solved with PENBMI for $N=3$ & $508$ & $1.05$\\

\rowcolor[gray]{.85}\centering Corollary \ref{corollary1} solved with PENBMI for $N=1$ & $425$ & $0.18$\\
\rowcolor[gray]{.85}\centering Corollary \ref{corollary1} solved with PENBMI for $N=2$ & $791$ & $0.23$\\
\rowcolor[gray]{.85}\centering Corollary \ref{corollary1} solved with PENBMI for $N=3$ & $925$ & $0.49$\\

\rowcolor[gray]{.85}\centering Algorithm 2 with $N=1$ & $427$ & $0.18$\\
\rowcolor[gray]{.85}\centering Algorithm 2 with $N=2$ & $642$ & $0.21$\\
\rowcolor[gray]{.85}\centering Algorithm 2 with $N=3$ & $842$ & $0.3$\\

\rowcolor[gray]{.85}\centering Algorithm 2 with $N=1$ and without constraint \eqref{corollary3-eq3} & $427$ & $0.17$\\
\rowcolor[gray]{.85}\centering Algorithm 2 with $N=2$ and without constraint \eqref{corollary3-eq3} & $673$ & $0.20$\\
\rowcolor[gray]{.85}\centering Algorithm 2 with $N=3$ and without constraint \eqref{corollary3-eq3} & $998$ & $0.28$\\

\rowcolor[gray]{.85}\centering Algorithm 3 with $N=1$ & $513$ & $5.46$\\
\rowcolor[gray]{.85}\centering Algorithm 3 with $N=2$ & $822$ & $8.24$\\
\rowcolor[gray]{.85}\centering Algorithm 3 with $N=3$ & $959$ & $12.47$\\

\rowcolor[gray]{.85}\centering Algorithm 3 with $N=1$ and without constraint \eqref{corollary3-eq3} & $513$ & $5.21$\\
\rowcolor[gray]{.85}\centering Algorithm 3 with $N=2$ and without constraint \eqref{corollary3-eq3} & $825$ & $9.86$\\
\rowcolor[gray]{.85}\centering Algorithm 3 with $N=3$ and without constraint \eqref{corollary3-eq3} & $999$ & $6.72$\\

\hline
\end{tabularx}
\label{table1}
\end{center}
\end{table*}
\end{example}

\begin{example}\label{example2}
In this example, we consider the discrete-time two-mass-spring
system from \cite{Kothare1996} with $(A,\,B,\,C^T ) = ({\bf{I}}_n
+ T_s A_c ,\,T_s B_c ,\,C_c^T )$, where
\begin{align*}
&(A_c ,\,B_c ,\,C_c^T ): = \left( {\left[ {\begin{array}{*{20}c}
   0 & 0 & 1 & 0  \\
   0 & 0 & 0 & 1  \\
   { - \frac{K}{{m_1 }}} & {\frac{K}{{m_1 }}} & 0 & 0  \\
   {\frac{K}{{m_2 }}} & { - \frac{K}{{m_2 }}} & 0 & 0  \\
\end{array}} \right],\,\left[ {\begin{array}{*{20}c}
   0  \\
   0  \\
   {\frac{1}{{m_1 }}}  \\
   0  \\
\end{array}} \right],\,\left[ {\begin{array}{*{20}c}
   1  \\
   0  \\
   0  \\
   1  \\
\end{array}} \right]} \right),
\end{align*}
$m_1$ and $m_2$ are the masses of the two bodies, $K$ is the
spring constant, and $T_s$ is the sampling time. The model
parameters are chosen to be $(m_1 ,\,m_2 ,\,K,\,T_s ) =
(1,\,1,\,1,\,0.05)$, and the open-loop system is unstable since
$\rho (A) = 1.0028$. For this system, the proposed two-steps and
ILMI algorithms with $N=1$ failed to find a solution. After
applying the proposed two-steps algorithm with $N=2$, a
$N$-PTVMSOFC gain matrix
\begin{align*}
&\hat {\cal F}_{{\rm{SOF}}}^{(2,\, \downarrow )}  = \left[
{\begin{array}{*{20}c}
   { - 167.7433} & 0  \\
   {460.2808} & { - 267.8199}  \\
\end{array}} \right]
\end{align*}
was obtained. Moreover, using Algorithm 1, the equivalent
closed-loop LTI system matrix was calculated to be
\begin{align*}
&{\cal A}_{{\rm{LTI}}} (\hat {\cal F}_{{\rm{SOF}}}^{(2,\,
\downarrow )} ,\,\Sigma ) = \left[ {\begin{array}{*{20}c}
   {0.5781} & {0.0025} & {0.1} & { - 0.4194}  \\
   {0.0025} & {0.9975} & 0 & {0.1}  \\
   {0.4663} & {0.7695} & {0.3280} & {1.2384}  \\
   {0.1} & { - 0.1} & {0.0025} & {0.9975}  \\
\end{array}} \right]
\end{align*}
with $\rho ({\cal A}_{{\rm{LTI}}} (\hat {\cal
F}_{{\rm{SOF}}}^{(2,\, \downarrow )} ,\,\Sigma )) = 0.9529$ and
eigenvalues $(0.2424,\,0.7602,\,0.9492 \pm 0.0754i)$. Finally, the
simulation result with $x(0) = [\begin{array}{*{20}c}
   3 & { - 3} & 0 & 0  \\ \end{array}]^T $ is depicted in Fig. \ref{figure1}.
\begin{figure}[t]
\centering 
%
%
\begin{psfrags}%
\psfragscanon%
%
\psfrag{s09}[t][t]{\color[rgb]{0,0,0}\setlength{\tabcolsep}{0pt}\begin{tabular}{c}$k$\end{tabular}}%
\psfrag{s10}[b][b]{\color[rgb]{0,0,0}\setlength{\tabcolsep}{0pt}\begin{tabular}{c}$x(k)$\end{tabular}}%
\psfrag{s14}[][]{\color[rgb]{0,0,0}\setlength{\tabcolsep}{0pt}\begin{tabular}{c} \end{tabular}}%
\psfrag{s15}[][]{\color[rgb]{0,0,0}\setlength{\tabcolsep}{0pt}\begin{tabular}{c} \end{tabular}}%
\psfrag{s16}[l][l]{\color[rgb]{0,0,0}$x_4(k)$}%
\psfrag{s17}[l][l]{\color[rgb]{0,0,0}$x_1(k)$}%
\psfrag{s18}[l][l]{\color[rgb]{0,0,0}$x_2(k)$}%
\psfrag{s19}[l][l]{\color[rgb]{0,0,0}$x_3(k)$}%
\psfrag{s20}[l][l]{\color[rgb]{0,0,0}$x_4(k)$}%
%
\psfrag{x01}[t][t]{$0$}%
\psfrag{x02}[t][t]{$0.1$}%
\psfrag{x03}[t][t]{$0.2$}%
\psfrag{x04}[t][t]{$0.3$}%
\psfrag{x05}[t][t]{$0.4$}%
\psfrag{x06}[t][t]{$0.5$}%
\psfrag{x07}[t][t]{$0.6$}%
\psfrag{x08}[t][t]{$0.7$}%
\psfrag{x09}[t][t]{$0.8$}%
\psfrag{x10}[t][t]{$0.9$}%
\psfrag{x11}[t][t]{$1$}%
\psfrag{x12}[t][t]{$0$}%
\psfrag{x13}[t][t]{$10$}%
\psfrag{x14}[t][t]{$20$}%
\psfrag{x15}[t][t]{$30$}%
\psfrag{x16}[t][t]{$40$}%
\psfrag{x17}[t][t]{$50$}%
\psfrag{x18}[t][t]{$60$}%
\psfrag{x19}[t][t]{$70$}%
\psfrag{x20}[t][t]{$80$}%
\psfrag{x21}[t][t]{$90$}%
\psfrag{x22}[t][t]{$100$}%
%
\psfrag{v01}[r][r]{$0$}%
\psfrag{v02}[r][r]{$0.1$}%
\psfrag{v03}[r][r]{$0.2$}%
\psfrag{v04}[r][r]{$0.3$}%
\psfrag{v05}[r][r]{$0.4$}%
\psfrag{v06}[r][r]{$0.5$}%
\psfrag{v07}[r][r]{$0.6$}%
\psfrag{v08}[r][r]{$0.7$}%
\psfrag{v09}[r][r]{$0.8$}%
\psfrag{v10}[r][r]{$0.9$}%
\psfrag{v11}[r][r]{$1$}%
\psfrag{v12}[r][r]{$-20$}%
\psfrag{v13}[r][r]{$-15$}%
\psfrag{v14}[r][r]{$-10$}%
\psfrag{v15}[r][r]{$-5$}%
\psfrag{v16}[r][r]{$0$}%
\psfrag{v17}[r][r]{$5$}%
\psfrag{v18}[r][r]{$10$}%
\psfrag{v19}[r][r]{$15$}%
\psfrag{v20}[r][r]{$20$}%
%
\epsfig{file=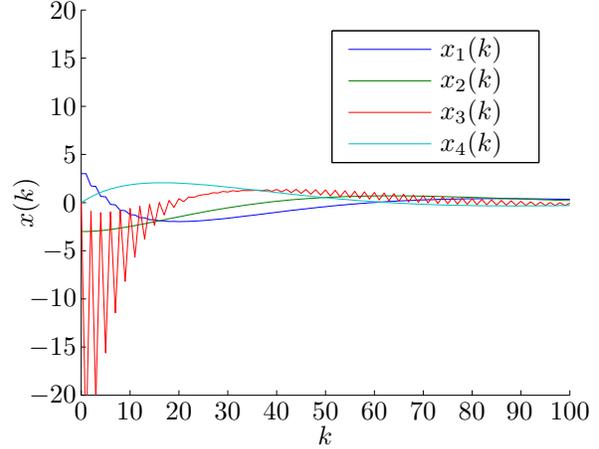,width=7.5cm}%
\end{psfrags}%
%
\caption{Example
\ref{example2}. Time histories of the state variables $x_1
(k),\,x_2 (k),\,x_3 (k)$, and $x_4 (k)$.}\label{figure1}
\end{figure}
\end{example}

\subsection{Reduction of chattering effect}

\begin{figure}[t]
\centering 
%
%
\begin{psfrags}%
\psfragscanon%
%
\psfrag{s09}[t][t]{\color[rgb]{0,0,0}\setlength{\tabcolsep}{0pt}\begin{tabular}{c}$k$\end{tabular}}%
\psfrag{s10}[b][b]{\color[rgb]{0,0,0}\setlength{\tabcolsep}{0pt}\begin{tabular}{c}$x(k)$\end{tabular}}%
\psfrag{s14}[][]{\color[rgb]{0,0,0}\setlength{\tabcolsep}{0pt}\begin{tabular}{c} \end{tabular}}%
\psfrag{s15}[][]{\color[rgb]{0,0,0}\setlength{\tabcolsep}{0pt}\begin{tabular}{c} \end{tabular}}%
\psfrag{s16}[l][l]{\color[rgb]{0,0,0}$x_4(k)$}%
\psfrag{s17}[l][l]{\color[rgb]{0,0,0}$x_1(k)$}%
\psfrag{s18}[l][l]{\color[rgb]{0,0,0}$x_2(k)$}%
\psfrag{s19}[l][l]{\color[rgb]{0,0,0}$x_3(k)$}%
\psfrag{s20}[l][l]{\color[rgb]{0,0,0}$x_4(k)$}%
%
\psfrag{x01}[t][t]{$0$}%
\psfrag{x02}[t][t]{$0.1$}%
\psfrag{x03}[t][t]{$0.2$}%
\psfrag{x04}[t][t]{$0.3$}%
\psfrag{x05}[t][t]{$0.4$}%
\psfrag{x06}[t][t]{$0.5$}%
\psfrag{x07}[t][t]{$0.6$}%
\psfrag{x08}[t][t]{$0.7$}%
\psfrag{x09}[t][t]{$0.8$}%
\psfrag{x10}[t][t]{$0.9$}%
\psfrag{x11}[t][t]{$1$}%
\psfrag{x12}[t][t]{$0$}%
\psfrag{x13}[t][t]{$10$}%
\psfrag{x14}[t][t]{$20$}%
\psfrag{x15}[t][t]{$30$}%
\psfrag{x16}[t][t]{$40$}%
\psfrag{x17}[t][t]{$50$}%
\psfrag{x18}[t][t]{$60$}%
\psfrag{x19}[t][t]{$70$}%
\psfrag{x20}[t][t]{$80$}%
\psfrag{x21}[t][t]{$90$}%
\psfrag{x22}[t][t]{$100$}%
%
\psfrag{v01}[r][r]{$0$}%
\psfrag{v02}[r][r]{$0.1$}%
\psfrag{v03}[r][r]{$0.2$}%
\psfrag{v04}[r][r]{$0.3$}%
\psfrag{v05}[r][r]{$0.4$}%
\psfrag{v06}[r][r]{$0.5$}%
\psfrag{v07}[r][r]{$0.6$}%
\psfrag{v08}[r][r]{$0.7$}%
\psfrag{v09}[r][r]{$0.8$}%
\psfrag{v10}[r][r]{$0.9$}%
\psfrag{v11}[r][r]{$1$}%
\psfrag{v12}[r][r]{$-20$}%
\psfrag{v13}[r][r]{$-15$}%
\psfrag{v14}[r][r]{$-10$}%
\psfrag{v15}[r][r]{$-5$}%
\psfrag{v16}[r][r]{$0$}%
\psfrag{v17}[r][r]{$5$}%
\psfrag{v18}[r][r]{$10$}%
\psfrag{v19}[r][r]{$15$}%
\psfrag{v20}[r][r]{$20$}%
%
\epsfig{file=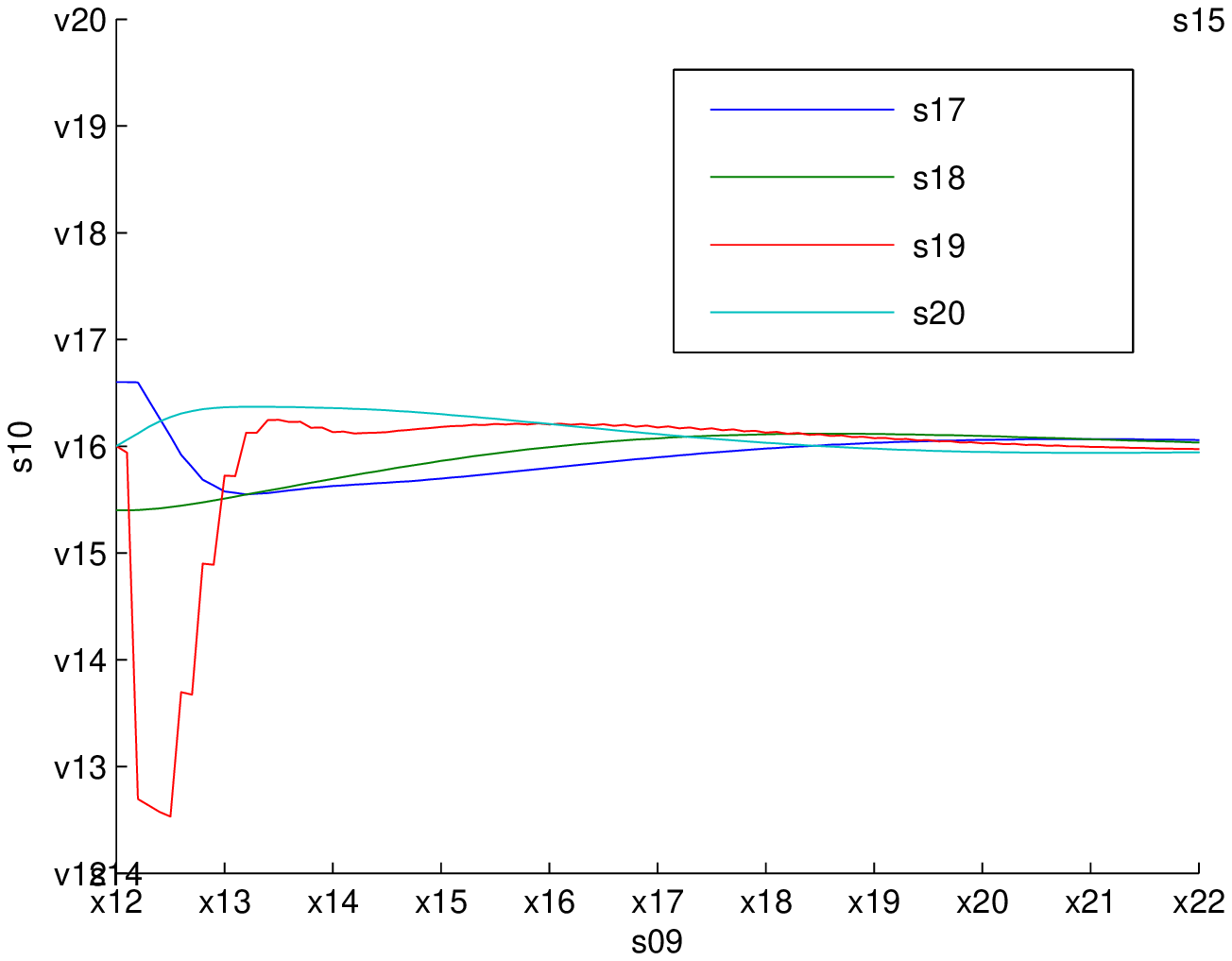,width=7.5cm}%
\end{psfrags}%
%
\caption{Example
\ref{example3}. Time histories of the state variables $x_1
(k),\,x_2 (k),\,x_3 (k)$, and $x_4 (k)$.}\label{figure2}
\end{figure}

As we can observe from Table \ref{table1}, the proposed method
outperforms the existing SOF approaches. Unfortunately, typically
the performance of the $N$-PTVMSOFC might not be so good since the
asymptotic stability is guaranteed only for the states
$x(k),\,\forall k \in \{ k \in {\mathbb N}:\,\left\lceil k
\right\rceil _N  = 0\}$. This property can cause chattering
problems as we can see from Fig. \ref{figure1}. In order to
alleviate the problem, we propose a simple procedure which may be
helpful to some degree in reducing the chattering effect.
Specifically, once a solution to the $N$-PTVMSOFC design is
obtained, then $\rho ({\cal A}_{{\rm{LTI}}} ({\cal
F}_{{\rm{SOF}}}^{(N,\, \downarrow )} ,\,\Sigma )) < 1$ holds. On
the other hand, the constraint is not satisfied for ${\cal
A}_{{\rm{LTI}}} ({\cal F}_{{\rm{SOF}}}^{(i ,\, \downarrow )}
,\,\Sigma ),\,i \in {\mathbb Z}_{[1,\,N - 1]}$ that correspond to
the intermediate states between $x(k)$ and $x(k+N)$. If $\rho
({\cal A}_{{\rm{LTI}}} ({\cal F}_{{\rm{SOF}}}^{(i ,\, \downarrow
)} ,\,\Sigma )) < \beta ,\,i \in {\mathbb Z}_{[1,\,N - 1]}$, then
a larger $\beta$ means larger fluctuations of states between
$x(k)$ and $x(k+N)$. In this perspective, we can try to minimize
$\beta$ while imposing constraint $\rho ({\cal A}_{{\rm{LTI}}}
({\cal F}_{{\rm{SOF}}}^{(N,\, \downarrow )} ,\,\Sigma )) < 1$. In
the sense of Lyapunov, this problem is equivalent to minimizing
$\beta$ subject to $P \succ 0,\,S \succ 0$, and
\begin{align*}
&{\cal A}_{{\rm{LTI}}} ({\cal F}_{{\rm{SOF}}}^{(N,\, \downarrow )}
,\,\Sigma )^T P{\cal A}_{{\rm{LTI}}} ({\cal F}_{{\rm{SOF}}}^{(N,\,
\downarrow )} ,\,\Sigma ) - P \prec 0,\\
&{\cal A}_{{\rm{LTI}}} ({\cal F}_{{\rm{SOF}}}^{(i ,\, \downarrow
)} ,\,\Sigma )^T S{\cal A}_{{\rm{LTI}}} ({\cal F}_{{\rm{SOF}}}^{(i
,\, \downarrow )} ,\,\Sigma ) - \beta S
\prec 0,\\
&i \in {\mathbb Z}_{[1,\,N - 1]}.
\end{align*}
Based on Corollary \ref{corollary3}, we can readily arrive at the
following result, which is presented without the proof:
\begin{corollary}\label{corollary4}
Suppose that $(\hat P,\,\hat {\cal J},\,\hat {\cal G})$ is a
solution to \eqref{Ebihara2011-lemma-eq1}, and let $\hat {\cal
F}_{{\rm{SF}}}^{(N,\, \downarrow )}  = ({\bf{T}}_N \otimes
{\bf{I}}_m )\hat {\cal J}\hat {\cal G}^{- 1} ({\bf{T}}_N \otimes
{\bf{I}}_n )$. Then, system \eqref{system} is stabilizable via
$N$-PTVMSOFC \eqref{PTVMSOFC} and $\rho ({\cal A}_{{\rm{LTI}}}
({\cal F}_{{\rm{SOF}}}^{(i ,\, \downarrow )} ,\,\Sigma )) < \beta
,\,i \in {\mathbb Z}_{[1,\,N - 1]}$ are guaranteed if there exists
matrices $P = P^T  \in {\mathbb R}^{n \times n} ,\,S = S^T \in
{\mathbb R}^{n \times n} ,\,M^{(i,\,j)}  \in {\mathbb R}^{m \times
p} ,\,V_{11}^{(i,\,j)} \in {\mathbb R}^{m \times m} ,\,{\cal
V}_{12}  \in {\mathbb R}^{Nm \times (N - 1)n}$, and ${\cal
V}_{22}\in {\mathbb R}^{(N - 1)n \times (N - 1)n}$ such that
\eqref{corollary3-eq1}, \eqref{corollary3-eq2},
\eqref{corollary3-eq3} with $\gamma = 1$, and the following LMI
problem is satisfied:
\begin{align}
&S \succ 0,\label{corollary4-eq1}\\
&\Pi _N^T {\cal X}_i (S,\,\beta )\Pi _N+ {\rm{He}}\{ {\cal H}(\hat
{\cal F}_{{\rm{SF}}}^{(N,\, \downarrow )} )^T {\cal D}({\cal
M},\,{\cal V}_{11} ,\,{\cal V}_{12} ,\,{\cal V}_{22}
)\}\nonumber\\
&\quad \prec 0,\quad i  \in {\mathbb Z}_{[1,\,N -
1]},\label{corollary4-eq2}
\end{align}
where ${\cal V}_{11}\in {\mathbb R}^{Nm \times Nm}$ and ${\cal
M}\in {\mathbb R}^{Nm \times Np}$ are defined in \eqref{V11} and
\eqref{M}, respectively. Moreover, an admissible PTVMSOFC gain
matrix is given by ${\cal F}_{{\rm{SOF}}}^{(N,\, \downarrow )}  =
{\cal V}_{11}^{- 1} {\cal M}$.
\end{corollary}
To reduce $\beta$, an ILMI algorithm similar to Algorithm 3 can be
applied based on Corollary \ref{corollary4}, although it is not
addressed here for space limitations.
\begin{example}\label{example3}
Let us consider Example \ref{example2} again. For $\hat {\cal
F}_{{\rm{SOF}}}^{(2,\, \downarrow )}$ given in Example
\ref{example2}, we applied an ILMI algorithm to reduce $\beta$ and
obtained gain matrix
\begin{align*}
&\hat {\cal F}_{{\rm{SOF}}}^{(2,\, \downarrow )}  = \left[
{\begin{array}{*{20}c}
   {0.0018} & 0  \\
   {365.0515} & { - 428.3888}  \\
\end{array}} \right]
\end{align*}
with $\rho ({\cal A}_{{\rm{LTI}}} (\hat {\cal
F}_{{\rm{SOF}}}^{(2,\, \downarrow )} ,\,\Sigma )) = 0.9535$. The
eigenvalues of ${\cal A}_{{\rm{LTI}}} (\hat {\cal
F}_{{\rm{SOF}}}^{(2,\, \downarrow )} ,\,\Sigma )$ were $(0.5094
\pm 0.3349i,\,0.9501 \pm 0.0806i)$ and $\rho ({\cal
A}_{{\rm{LTI}}} (\hat {\cal F}_{{\rm{SOF}}}^{(1,\, \downarrow )}
,\,\Sigma )) = 1.0025$, while in Example \ref{example2}, $\rho
({\cal A}_{{\rm{LTI}}} (\hat {\cal F}_{{\rm{SOF}}}^{(1,\,
\downarrow )} ,\,\Sigma )) = 1.2141$. The simulation result under
the same initial condition is plotted in Fig. \ref{figure2}, which
clearly shows that the amplitude of oscillation was mitigated in
comparison with that of Fig. \ref{figure1}.
\end{example}
\begin{remark}
An extension of our methods to the LQR formulation is
straightforward. Let us consider the following cost function:
\begin{align*}
J_\infty  (x,\,u):=& \sum\limits_{k \in \{ k \in {\mathbb
N}:\,\left\lceil k \right\rceil _N  = 0\} } {\left[
{\begin{array}{*{20}c}
   {x(k:k + N - 1)}  \\
   {u(k:k + N - 1)}  \\
\end{array}} \right]^T }\\
&\times W\left[ {\begin{array}{*{20}c}
   {x(k:k + N - 1)}  \\
   {u(k:k + N - 1)}  \\
\end{array}} \right],
\end{align*}
where $W: = \left[ {\begin{array}{*{20}c}
   Q & {\bf{0}}  \\
   {\bf{0}} & R  \\
\end{array}} \right] \succeq 0$ is a given weighting
matrix. Then, with only a little modification, it is easy to see
that if the LMIs of Corollary \ref{corollary3} with
\eqref{corollary3-eq2} replaced by
\begin{align*}
&\Pi _N^T {\cal X}_N (P,\,1)\Pi _N  + W\\
&\quad + {\rm{He}}\{ {\cal H}(\hat {\cal F}_{{\rm{SF}}}^{(N,\,
\downarrow )} )^T {\cal D}({\cal M},\,{\cal V}_{11} ,\,{\cal
V}_{12} ,\,{\cal V}_{22} )\}  \prec 0,
\end{align*}
is satisfied, then $N$-periodic control system
\eqref{SOF-closed-loop-system} with ${\cal F}_{{\rm{SOF}}}^{(N,\,
\downarrow )}  = {\cal V}_{11}^{ - 1} {\cal M}$ is asymptotically
stable, and the cost function satisfies the bound $J_\infty
(x,\,u) < x(0)^T Px(0)$.
\end{remark}


\end{document}